\documentclass{article}
\usepackage{latexsym,amsfonts,amsmath,amsthm,amssymb,makeidx}
\usepackage[title]{appendix}
\usepackage{CJK,CJKnumb,CJKulem,times,dsfont,ifthen,mathrsfs,latexsym,amsfonts, color}
\usepackage{amsmath,amsthm,makeidx,fontenc,amssymb,bm,graphicx,psfrag,listings, curves,extarrows}

\let\oldbibliography\thebibliography
\renewcommand{\thebibliography}[1]{%
\oldbibliography{#1}%
\setlength{\itemsep}{0pt}%
}

\makeindex
\newtheorem{definition}{Definition}[section]
\newtheorem{theorem}{Theorem}[section]
\newtheorem{lemma}{Lemma}[section]
\newtheorem{corollary}{Corollary}[section]

\newtheorem{remark}{Remark}[section]

\newcommand{\q}{\theta}

\newcommand{\s}{\section}

\newcommand{\la}{\lambda}

\newcommand{\pa}{\partial}

\newcommand{\R}{\mathbb R}

\newcommand{\bt}{\begin{theorem}}
\newcommand{\et}{\end{theorem}}
\newcommand{\bl}{\begin{lemma}}
\newcommand{\el}{\end{lemma}}
\newcommand{\bd}{\begin{definition}}
\newcommand{\ed}{\end{definition}}
\newcommand{\bc}{\begin{corollary}}
\newcommand{\ec}{\end{corollary}}
\newcommand{\bp}{\begin{proof}}
\newcommand{\ep}{\end{proof}}
\newcommand{\bx}{\begin{example}}
\newcommand{\ex}{\end{example}}
\newcommand{\bi}{\begin{exercise}}
\newcommand{\ei}{\end{exercise}}
\newcommand{\bo}{\begin{prop}}
\newcommand{\eo}{\end{prop}}
\newcommand{\br}{\begin{remark}}
\newcommand{\er}{\end{remark}}
\newcommand{\be}{\begin{equation}}
\newcommand{\ee}{\end{equation}}
\newcommand{\ba}{\begin{align}}
\newcommand{\ea}{\end{align}}
\newcommand{\bn}{\begin{enumerate}}
\newcommand{\en}{\end{enumerate}}
\newcommand{\bg}{\begin{align*}}
\newcommand{\bcs}{\begin{cases}}
\newcommand{\ecs}{\end{cases}}

\newcommand{\bean}{\begin{eqnarray*}}
\newcommand{\eean}{\end{eqnarray*}}

\numberwithin{equation}{section}

\begin{document}

\title{{\bf On the  Triharmonic  Lane-Emden Equation}\thanks{Supported by NSFC of China and NSERC of Canada}}

\date{}
\author{\\{\bf  Senping Luo$^{1}$,\;\;  Juncheng Wei$^{2}$ \;\; and \; Wenming Zou$^{3}$}\\
\footnotesize {\it  $^{1,3}$Department of Mathematical Sciences, Tsinghua University, Beijing 100084, China}\\
\footnotesize {\it  $^{2}$Department of Mathematics, University of British Columbia, Vancouver, BC V6T 1Z2,
Canada}
}

\maketitle
\begin{center}
\begin{minipage}{120mm}
\begin{center}{\bf Abstract}\end{center}

We derive a monotonicity formula and classify  finite Morse index solutions (positive or sign-changing, radial or not) to  the following  triharmonic Lane-Emden equation:
\noindent
\begin{equation}\nonumber
(-\Delta)^3 u=|u|^{p-1}u\;\;\;\;\hbox{in}\;\;\;\;\; \R^n,
\end{equation}
where $p$ is below the Joseph-Lundgren exponent. As a byproduct we also obtain a new monotonicity formula for the triharmonic maps.

\vskip0.10in


\end{minipage}

\end{center}

\vskip0.10in
\section{Introduction and Main results}

In this paper,  we study the finite Morse index solutions of the following triharmonic Lane-Emden equation
 \begin{equation}\label{10LE}
(-\Delta)^3 u=|u|^{p-1}u\;\;\;\;\hbox{in}\;\;\;\;\;  \R^n
\end{equation}
 and give a complete classification of such kind of solutions.

The   Lane-Emden equation
 \begin{equation}\label{LZ=0001}
 -\Delta  u=|u|^{p-1}u \;\;\;\;\hbox{in}\;\;\;\;\;  \R^n
\end{equation}
and its parabolic counterpart have played an essential role in the development of methods of nonlinear PDEs in the last
decades. A fundamental result on equation \eqref{LZ=0001} is the celebrated   Liouville-type theorem due to Gidas and Spruck \cite{Gidas-Sp=1981}: {\it The     equation \eqref{LZ=0001} has no positive  classical solution if $0<p <p_S$, where $p_S :=(n+2)/(n-2)$ if $ n\geq 3$ while  $p_S :=\infty$ if $n\leq 2.$ } Since then  there has been  an extensive literature on such a  type of equations  or systems. In particular, in 2007 the seminar paper \cite{Farina2007} by Farina  (see also \cite{Farina2005}), the equation  \eqref{LZ=0001} is revisited  for $p>\frac{n+2}{n-2}$.  The author obtained  some  classification results and Liouville-type theorems  for smooth solutions  including stable solutions, finite Morse index solutions, solutions which are stable
outside a compact set, radial solutions and non-negative solutions. The  results obtained in  \cite{Farina2007} were   applied to subcritical, critical and supercritical values of the exponent $p.$  Moreover,   the critical stability exponent $p_c(n)$ (Joseph-Lundgren exponent) is determined which   is larger than the classical critical exponent  $p_S=2^\ast-1$ in Sobolev imbedding theorems. In some sense, the   Joseph-Lundgren exponent  $p_c(n)$ is a critical threshold for obtaing the  Liouville-type theorems for stable or finite Morse index solutions.  The proof of Farina involves a delicate use of Nash-Moser's iteration technique, which is a classical tool for regularity of second order elliptic operators and falls short for higher order operators.

The  biharmonic  Lane-Emden equation:
\begin{equation}
  (-\Delta)^2 u=|u|^{p-1}u\quad \hbox{ in }\;   \R^n
  \end{equation}
  has also attracted lots of studies in recent years. The classical Gidas-Spruck type result has been extended (\cite{Lin1998}, \cite{Wei1999}). The radial solutions are classified (\cite{Gazzola2006}, \cite{GW2010}). The classification of stable/finite Morse index solutions was initiated by Cowan-Esposito-Ghoussoub  \cite{CEG} and Cowan-Ghoussoub \cite{CG}. A complete classification was obtained by  Davila-Dupaigne-Wang-Wei in \cite{Wei=2}. They
give a complete classification of stable and finite Morse index solutions (whether positive or sign changing), in the full exponent range.  To by-pass the Nash-Moser iterations, a key point used in \cite{Wei=2} is the monotonicity formula for bi-harmonic equations.

   On the other hand,  very recently  the  nonlocal  Lane-Emden equation
 \begin{equation}\label{10LE=zou}
(-\Delta)^s u=|u|^{p-1}u\;\;\;\;\hbox{ in }  \R^n
\end{equation}
were considered in Davila-Dupaigne-Wei  \cite{Wei0=1} when $0<s<1$ and Fazly-Wei   \cite{Wei1=2} when $1<s<2$. Both \cite{Wei0=1} and \cite{Wei1=2} gave a complete classification of finite Morse index solution of \eqref{10LE=zou}.


The motivation to study the above equations   comes from   both
physics and   geometry. In particular, the critical case is  inevitable for studying the  conformal geometry like  the prescribed
scalar curvature problem. On the other hand, it is well known that the Liouville-type theorems  play a crucial role to get  a priori $L^\infty$-bounds for
solutions of semilinear   elliptic and parabolic problems. In this regard,  we refer to the book by Quittner and Souplet \cite{QS2007}.

In this paper  we initiate  the  study of finite Morse index solutions to  the   triharmonic Lane-Emden equation  \eqref{10LE}. There are three critical exponents. The first one is the Serrin's exponent $ \frac{n}{n-6}$. The second is the Sobolev exponent $p_S=\frac{n+6}{n-6}$. The third is   the Joseph-Lundgren exponent which is  given by the following formula:
\be\label{10pcn}
p_c(n)=
\begin{cases}
\;\;\;\infty\;\;&\hbox{if}\;\; n\leq 14,\\
\frac{n+4-2d(n)}{n-8-2d(n)}\;\;&\hbox{if}\;\; n\geq15,
\end{cases}
\ee
where
\be\label{dndef}
d(n):=\frac{1}{6}\Big(9n^2+96-\frac{1536+1152n^2}{d_0(n)}-\frac{3}{2}d_0(n)\Big)^{1/2};
\ee
\be\nonumber
d_0(n):=-(d_1(n)+36\sqrt{d_2(n)})^{1/3};
\ee
\be\nonumber\aligned
d_1(n):=-94976+20736n+103104n^2-10368n^3+1296n^5-3024n^4-108n^6;
\endaligned\ee
\be\nonumber\aligned
d_2(n):&=6131712-16644096n^2+6915840n^4-690432n^6-3039232n\\
&\quad+4818944n^3-1936384n^5+251136n^7-30864n^8-4320n^9\\
&\quad+1800n^{10}-216n^{11}+9n^{12}.
\endaligned\ee

\br  In the harmonic case, the Joseph-Lundgren exponent  (Joseph-Lundgren \cite{Joseph1972})  is given by
\be
p_{cHarmonic}(n):=\begin{cases}
\;\;\;\;\;\;\;\;\;\infty\;\;\;\;\;\;\;\;&\hbox{if }\;  n\leq10,\\
\frac{(n-2)^2-4n+8\sqrt{n-1}}{(n-2)(n-10)}\;\;\;\;\;\;\;\;&\hbox{if }\;  n\geq11
\end{cases}
\ee
while in the bi-harmonic case, the corresponding exponent (Gazzola and Grunau \cite{Gazzola2006}) is
\be
p_{cBiharmonic}(n):=\begin{cases}
\;\;\;\;\;\;\;\;\;\;\infty\;\;\;\;\;\;\;\;&\hbox{if}\;  n\leq12,\\
\frac{n+2-\sqrt{n^2+4-n\sqrt{n^2-8n+32}}}{n-6-\sqrt{n^2+4-n\sqrt{n^2-8n+32}}}\;\;\;\;\;\;\;\;&\hbox{if}\;  n\geq13.
\end{cases}
\ee
 In the triharmonic case, $p_c (n)$ satisfies a $6$-th order polynomial algebraic  equation which in general has no explicit solution. It is interesting that we obtain explicit formula.
\er

\vskip0.1in

Next, we recall several definitions.

\bd\label{zou-stable} A solution $u$ of \eqref{10LE} is said to be stable if
\be\nonumber
\int_{\R^n}|\nabla\Delta \varphi|^2dx\geq p\int_{\R^n}|u|^{p-1}\varphi^2dx,\quad\hbox{for any}\;\varphi\in H^3 (\R^n).
\ee
\ed

\bd\label{zou-stable-out} A solution $u$ of \eqref{10LE} is said to be stable outside a compact $\Theta\subset \R^n $if
\be\nonumber
\int_{\R^n}|\nabla\Delta \varphi|^2dx\geq p\int_{\R^n}|u|^{p-1}\varphi^2dx,\quad\hbox{for any}\;\varphi\in H^3(\R^n\backslash\Theta).
\ee
\ed

\bd The Morse index   of the solution $u$ of \eqref{10LE} is defined as the maximal dimension over  all subspaces $E$ of $H^3 (\R^n)$ satisfying
\be\nonumber
\int_{\R^n}|\nabla\Delta \varphi|^2dx<p\int_{\R^n}|u|^{p-1}\varphi^2dx,\quad\hbox{for any}\; \varphi\in E\setminus\{0\}.
\ee
\ed
Hence, a solution is stable if and only if its Morse index is equal to zero. It is known that if a solution $u$ to \eqref{10LE} has finite Morse index, then there
exists a compact set $\mathcal{K}\subset\R^n$ such that
\be\nonumber
\int_{\R^n}|\nabla\Delta \varphi|^2dx\geq p\int_{\R^n}|u|^{p-1}\varphi^2dx,\quad\hbox{for any} \;\varphi\in H^2(\R^n\setminus \mathcal{K}).
\ee

The first main result  of the present paper is the following

\bt\label{10Liouville}
 Let $u$ be a stable solution of \eqref{10LE}.  If $ 1<p< p_c (n)$, then $ u \equiv 0$.

\et

For finite Morse index solutions we have the following

\bt\label{10Liouvillec}
Let $u$ be a finite Morse index solution of \eqref{10LE}.
Assume that either
\begin{itemize}
\item [(1)]   $1<p<\frac{n+6}{n-6}$
or
\item [(2)]  $\frac{n+6}{n-6}< p<p_c(n)$,
\end{itemize}
then the solution $u\equiv0$.
\begin{itemize}
\item[(3)]  If $p=\frac{n+6}{n-6}$, then $u$ has a finite energy, i.e.,
\be\nonumber
\int_{\R^n}|\nabla\Delta u|^2=\int_{\R^n}|u|^{p+1}<+\infty.
\ee

\end{itemize}
\et

\br
In both Theorems the condition $p <p_c (n)$ is optimal. In fact the radial singular solution is stable for $p\geq p_c (n)$. See \cite{LWZ}.
\er

\br While there are many works on second order and fourth order Lane-Emden equations, there are very few works on   $6-$th order Lane-Emden equations. We refer to  Farina-Ferreo \cite{Farina2016}, Lazzo-Schmidt \cite{LS2009},  Martinazzi \cite{Mar} and the references therein for related results on polyharmonic nonlinear equations.
\er

Theorems  \ref{10Liouville}-\ref{10Liouvillec} are proved by Monotonicity Formula which we introduce in the next section.

\s{Monotonicity formula for triharmonic Lane-Emden equations}

\vskip0.1in


We denote $\pa_r u= \nabla u\cdot\frac{x}{r}, r=|x|$.  Let $\delta_i,i=1,2,3,4$ be  defined by
\be\label{10delta}\aligned
\delta_1=&2n-\frac{24}{p-1},\\
\delta_2=&n(n-2)-n\frac{36}{p-1}-\frac{36}{p-1}(1+\frac{36}{p-1}),\\
\delta_3=&-\frac{24}{p-1}(1+\frac{6}{p-1})(2+\frac{6}{p-1})+2n\frac{12}{p-1}(1+\frac{6}{p-1})\\
           &-(n+b)(n+b-2)(1+\frac{12}{p-1}),\\
\delta_4=&(3+\frac{6}{p-1})(2+\frac{6}{p-1})(1+\frac{6}{p-1})\frac{6}{p-1}-2n(1+\frac{6}{p-1})(2+\frac{6}{p-1})\frac{6}{p-1}\\
          &+n(n-2)(2+\frac{6}{p-1})\frac{6}{p-1}.\\
\endaligned\ee

Next, we will introduce a functional and consider its monotonicity formula.  Let
$$B_\lambda:=\{y\in \R^n: |y-x|<\lambda\},\;\;  \lambda>0.$$
We define the functional $E(\la,x,u)$ depending on $x\in \R^n, \lambda>0$ and $ u$:
\be\label{10exu}\aligned
&E(\la,x,u)\\
&:=\la^{6\frac{p+1}{p-1}-n}\Big(\int_{ B_\la}\frac{1}{2}|\nabla\Delta u|^2
-\frac{1}{p+1}\int_{ B_\la}|u|^{p+1}\Big)\\
&\quad -\int_{\pa B_\la}\la^{6\frac{p+1}{p-1}-n-5}\Big[\frac{6}{p-1}(\frac{6}{p-1}-1)(\frac{6}{p-1}-2)(\frac{6}{p-1}-3)u\\
&\quad+\frac{24}{p-1}(\frac{6}{p-1}-1)(\frac{6}{p-1}-2)\la\pa_r u+\frac{36}{p-1}(\frac{6}{p-1}-1)\la^2\pa_{rr}u\\
&\quad+\frac{24}{p-1}\la^3\pa_{rrr}u+\la^4\pa_{rrrr}u\Big]\Big[\frac{6}{p-1}u+r\la\pa_r u\Big]\\
&\quad+2\int_{\pa B_\la}\la^{6\frac{p+1}{p-1}-n-5}\Big(\frac{6}{p-1}(\frac{6}{p-1}-1)(\frac{6}{p-1}-2)(\frac{6}{p-1}-3)u\\
&\quad+\frac{24}{p-1}(\frac{6}{p-1}-1)(\frac{6}{p-1}-2)\la\pa_r u+\frac{36}{p-1}(\frac{6}{p-1}-1)\la^2\pa_{rr}u\\
&\quad+\frac{24}{p-1}\la^3\pa_{rrr}u+\la^4\pa_{rrrr}u\Big)\Big(\frac{6}{p-1}u+r\la\pa_r u\Big)\\
&\quad-(\delta_1-6)\int_{\pa B_\la}\la^{6\frac{p+1}{p-1}-n-5}\big[\frac{6}{p-1}(\frac{6}{p-1}-1)(\frac{6}{p-1}-2)u\\
&\quad+\frac{18}{p-1}(\frac{6}{p-1}-1)\la\pa_r u
+\frac{18}{p-1}\la^2 \pa_{rr} u+\la^3\pa_{rrr} u\big]\\
&\quad\quad[\frac{6}{p-1}u+\la\pa_r u]+(\frac{6}{p-1}+2)\la^{6\frac{p+1}{p-1}-n-1}\int_{\pa B_\la}(\Delta_b u)^2\\
&\quad-(24-6\delta_1+\delta_2)\int_{\pa B_\la}\la^{6\frac{p+1}{p-1}-n-5}
\big[\frac{6}{p-1}(\frac{6}{p-1}-1) u\\
&\quad+\frac{14}{p-1}\la\pa_r u
+\la^2\pa_{rr}u\big]\big[\frac{6}{p-1}u+\la\pa_r u\big]\\
&\quad-(9\delta_1-3\delta_2-36)\int_{\pa B_\la}\la^{6\frac{p+1}{p-1}-n-5}
\big[\frac{6}{p-1}u+\la\pa_r u\big]^2\\
&\quad+(\delta_1-8)\int_{\pa B_\la}\la^{6\frac{p+1}{p-1}-n-5}
\big[\frac{6}{p-1}(\frac{6}{p-1}-1) u\\
&\quad+\frac{12}{p-1}\la\pa_r u+\la^2\pa_{rr}u\big]^2\\
&\quad+\delta_4\int_{\pa B_\la}\la^{6\frac{p+1}{p-1}-n-5}
\Big(\frac{6}{p-1}u+\la\pa_r u)\Big)u
\endaligned\ee


\be\nonumber\aligned
&+2\delta_4\int_{\pa B_\la}\la^{2s\frac{p+1}{p-1}-n-5}
u^2\\
&+2\int_{\pa B_\la}\la^{6\frac{p+1}{p-1}-n-5}(\la^2\Delta u-\la^2\pa_{rr}u-(n-1)\la\pa_r u)^2\\
&+\frac{1}{2}\int_{\pa B_\la}\la^{6\frac{p+1}{p-1}-n-4}\frac{d}{d\la}(\la^2\Delta u-\la^2\pa_{rr}u-(n-1)\la\pa_r u)^2\\
&-4(\beta-n+3)\int_{\pa B_\la}\la^{6\frac{p+1}{p-1}-n-5}\big[\la^2|\nabla u|^2-\la^2|\pa_r u|^2\big]\\
&-(\beta-n+3)\int_{\pa B_\la}\la^{6\frac{p+1}{p-1}-n-4}\frac{d}{d\la}\big[\la^2|\nabla u|^2-\la^2|\pa_r u|^2\big]\\
&-6\int_{\pa B_\la}\la^{6\frac{p+1}{p-1}-n-4}\frac{d}{d\la}\big[
\la^2|\nabla(\frac{6}{p-1}u+\la\pa_r u)|^2-\la^2|\pa_r(\frac{6}{p-1}u+\la\pa_r u)|^2
\big]^2\\
&+\int_{\pa B_\la}\la^{6\frac{p+1}{p-1}-n-2}\frac{d}{d\la}\big[
|\nabla(\frac{6}{p-1}u+\la\pa_r u)|^2-|\pa_r(\frac{6}{p-1}u+\la\pa_r u)|^2
\big]^2,\\
\endaligned\ee
where $\delta_j (j=1,2,3,4)$ are defined in \eqref{10delta},

\vskip0.1in

The following is the monotonicity formula which will play an important role.
\bt\label{10monoid}
Let $u$ satisfy  the equation \eqref{10LE}. Define $u^\la(x)=\la^{\frac{6}{p-1}}u(\la x)$,  then
\be\aligned
&\frac{d E(\la,x,u)}{d\la}\\
&=\int_{\pa B_1}\Big(2\la^5(\frac{d^3 u^\la}{d\la^3}\big)^2+(10\delta_1-2\delta_2-56)\la^3\big(\frac{d^2 u^\la}{d\la^2}\big)^2\\
&\quad +(-18\delta_1+6\delta_2-4\delta_3+2\delta_4+72)\la\big(\frac{d u^\la}{d\la}\big)^2\Big)\\
&\quad +\int_{\pa B_1}\Big(4\la^3\big(\frac{d^2}{d\la^2}\nabla_{\theta} u^\la\big)^2+(8\alpha-4\beta+4n-28)\la\big(\frac{d}{d\la}\nabla_{\theta} u^\la\big)^2\Big)\\
&\quad +2\la\int_{\pa B_1}\big[\frac{d}{d\la}\mathbf{div}_{\theta}(\nabla_{\theta}u^\la)\big]^2
+\la\int_{\pa B_1}\big(\frac{d\Delta u^\la}{d\la}\big)^2.\\
\endaligned
\ee
\et
We will give the proof of Theorem   \ref{10monoid} in the next section. Now we would like to
state a consequence  of  Theorem   \ref{10monoid}.

\vskip0.1in

 The functional $E(\la,x,u)$, defined in \eqref{10exu},  can be divided into two parts: the integral  over the ball $B_\la$
and the terms of integrals on the boundary $\pa B_\la$. We notice that in the blow-down analysis process,  the boundary terms
can be controlled  using initial  energy estimates. Then,  we may change some coefficients of the boundary terms in $E(\la,x,u)$, we denote it by $E^c(\la,x,u)$, which may be formulated  in the following way:
\be\label{LZ=623}\aligned&E^c(\la,x,u)\\
&:=E(\la,x,u)-\int_{\pa B_1}
\Big(\sum_{0\leq i,j\leq2}c_{i,j}\la^{i+j}\frac{d^iu^\la}{d\la^i}\frac{d^ju^\la}{d\la^j}
\Big),\endaligned\ee
where  $c_{i,j}\in \R$ are chosen  properly which may be different by various cases. Moreover, we still can obtain the lower bound of $\frac{d E^c(\la,x,u)}{d\la}$ and the lower bound
is independent of $c_{i,j}\in \R$. We have the following precise statement

\bt\label{10Monotone} Assume that $\frac{n+6}{n-6}<p<p_m (n) $. Then there exist $c_{ij}$ such that  $E^c(\la,x,u)$,  defined at (\ref{LZ=623}), is a nondecreasing function of $\la>0$. Furthermore,
\be\label{11mono}
\frac{d E^c(\la,x,u)}{d\la}\geq C(n,p)\la^{6\frac{p+1}{p-1}-6-n}
\int_{\pa B_\la(x_0)}\Big(\frac{6}{p-1}u+\la\pa_r u\Big)^2,
\ee
where $C(n,p)>0$ is a constant independent of $\la$, and
\be\nonumber
p_m(n):=\begin{cases}
+\infty\;\;\;\;\;\;\;\;&\hbox{if}\;\;\;\;\;\;\;\; n\leq30,\\
\frac{5n+30-\sqrt{15n^2-60n+190}}{5n-30-\sqrt{15n^2-60n+190}}&\hbox{if}\;\;\;\;\;\;\;\;n\geq31.\\
\end{cases}
\ee

\et

\br In the above theorem, we  need the upper bound condition of $p$, namely $p<p_m (n)$.  Let us recall that in the biharmonic case the monotonicity formula holds for  all $ p>\frac{n+4}{n-4}$ (see \cite{Wei=2}). Since $p_c(n)<p_m(n)$, the above monotonicity formula holds for $\frac{n+6}{n-6}<p<p_c (n) $ which is used for our blow down analysis. See Theorem \ref{LZ=62609}. It seems that in the triharmonic case, the supercritical condition  $p>\frac{n+6}{n-6}$ alone  is not sufficient to make such kind of monotonicity formula \eqref{11mono} hold. We refer the readers to  section 7 and \cite{Blatt} for more details.
\er

\br The proof of Theorem \ref{10Monotone} is quite involved. In the bi-harmonc cases, the positivity of $ \frac{dE}{d \lambda}$ is trivial. Here we have to discuss three cases: $n\leq20$, $21\leq n\leq30$ and $n\geq31$. In each case we have to come up with different combinations of terms.
\er

\br\label{10Blatt1}
In \cite{Blatt}, Simon Blatt also derived a monotonicity formula for triharmonic Lane-Emden equations under different conditions on $p$ which is much stronger than our's here (see below). He then used to prove partial regularity of stationary solutions  and obtain Hausdorff dimension estimates   for the singular set of solutions. By unifying  the notations, we see that  in
\cite{Blatt} (Corollary 3.13), the author gets the monotonicity under the condition $-(n-20)+8\alpha(n-1)\geq0$.  Transfer to our notations, the monotonicity formula of \cite{Blatt} requires  the hypothesis  $\frac{n+6}{n-6}<p<p_{m_1}(n)$, where
\be\nonumber
p_{m_1}(n):=\begin{cases}
+\infty\;\;\;\;\;\;\;\;&\hbox{if}\;\;\;\;\;\;\;\; n\leq20,\\
\frac{n+28}{n-20}&\hbox{if}\;\;\;\;\;\;\;\;n\geq21.\\
\end{cases}
\ee
A direct calculation shows that $p_{m_1}(n)<p_m(n)$. Therefore,  by using our arguments in the current paper,  the second main result  Theorem 1.2 of \cite{Blatt} actually  can be improved to
$\frac{n+6}{n-6}<p<p_{m}(n)$.

\er

By slightly modifying the proof of Theorem \ref{10Monotone}, we are able to get the monotonicity formula for the  triharmonic map, i.e.,
\be\nonumber
\Delta^3 u=0.
\ee
Indeed, let $p\rightarrow+\infty$ in \eqref{10exu} and denote $E_\infty(\la,x,u)=\lim_{p\rightarrow\infty}E(\la,x,u)$,   where the term $\frac{1}{p+1}\la^{6\frac{p+1}{p-1}-n}\int_{\pa B_\la}|u^\la|^{p+1}$ is understood   vanished, then we have
\bc
Assume that $7\leq n\leq30$. Then there exist $c_{ij}$ such that  $E^c_\infty(\la,x,u)$,  defined similarly  as in (\ref{LZ=623}), is a nondecreasing function of $\la>0$. Furthermore,
\be\label{11mono}
\frac{d E_\infty^c(\la,x,u)}{d\la}\geq C(n)\la^{-n}
\int_{\pa B_\la(x_0)}\Big(\la\pa_r u\Big)^2,
\ee
where $C(n)>0$ is a constant independent of $\la$.

\ec

\br\label{10Blatt2}
In \cite{Blatt}, Simon Blatt   derived a monotonicity formula for  extrinsic triharmonic maps under the conditions on $ 6< n\leq20$ by which a smoothness result was obtained (See Theorem 1.1 of \cite{Blatt}). However, by our results here, his Theorem 1.1 can be improved to $ 6< n\leq 30$.
\er

At the end of this section, we say a few words on the powerful applications of  monotonicity formula. It is known that
monotonicity formulas are one of the most important tools for studying geometric problems as well as supercritical equations
and systems. For monotonicity formulas for stationary harmonic maps we refer to Evans \cite{Evans} for harmonic maps and  Chang-Wang-Yang \cite{CWY} for biharmonic maps.    For the second order  Lane-Emden equation we refer to Giga-Kohn \cite{GK} and Pacard \cite{Pac}. For biharmonic and fractional Lane-Emden equations we refer to \cite{Wei0=1, Wei=2, Wei1=2}.


 \section{Monotonicity formula and the proof of Theorem   \ref{10monoid}}

Since the derivation of the derivative for the $E(\la,x,u)$ is complicated,  we divide it into several subsections.
In subsection $3.1$, we  derive  $\frac{d}{d\la}\overline{E}(u,\la)$, where $\overline{E}(u,\la)$ is defined in \eqref{zzz=9} below,
which in fact is the first term of $E(\lambda, x, u)$ introduced in  \eqref{10exu}.
In subsection $3.2$, we calculate the (higher-order) derivatives  $\frac{\pa^j}{\pa r^j}u^\la$  and  $\frac{\pa^i}{\pa \la^i}u^\la,\;\; i,j=1,2,3,4$.
In subsection $3.3$,  the operator   $\Delta^2$   and its representation will be given.
In subsection $3.4$, we  decompose $\frac{d}{d\la}\overline{E}(u^\la,1)$.
Finally, combining with the above four subsections, we can obtain the  derivative formula, hence get the proof of  Theorem \ref{10monoid}.

\vskip0.1in
Without loss of generality, suppose that $x_0=0$ and denote by $B_\la$  the ball   centered at zero with radius $\la$. Set
\be\label{zzz=9}
\overline{E}(u,\la):=\la^{6\frac{p+1}{p-1}-n}\Big(\int_{ B_\la}\frac{1}{2}|\nabla\Delta u|^2
-\frac{1}{p+1}\int_{ B_\la}|u|^{p+1}\Big).\\
\ee

\subsection{The derivation of $\frac{d}{d\la}\overline{E}(u,\la)$}
Define
\be\label{10uvw}\aligned
&v:=\Delta u, \quad\quad  u^\la(x):=\la^{\frac{6}{p-1}}u(\la x), \quad w:=\Delta v,\\
&v^\la(x):=\la^{\frac{6}{p-1}+2}v(\la x),\quad\quad   w^\la(x):=\la^{\frac{6}{p-1}+4}w(\la x).
\endaligned\ee
Therefore,
\be\label{10lamda}
\Delta u^\la(x)=v^\la(x), \quad\quad  \Delta v^\la(x)=w^\la(x).
\ee
In addition, differentiating \eqref{10lamda} with respect to $\la$ we have
\be\nonumber
\Delta \frac{du^{\la}}{d\la}=\frac{dv^{\la}}{d\la},\;\;\; \Delta \frac{dv^{\la}}{d\la}=\frac{dw^{\la}}{d\la}.
\ee
Note that
\be\nonumber
\overline{E}(u,\la)=\overline{E}(u^\la,1)=\int_{ B_1}\frac{1}{2}|\nabla v^\la|^2
-\frac{1}{p+1}\int_{B_1}|u^\la|^{p+1}.
\ee
Taking derivative of the energy $\overline{E}(u^\la,1)$ with respect to $\la$ and integrating by part, we have :

\be\label{10overE}\aligned
\frac{d\overline{E}(u^\la,1)}{d\la}&=\int_{ B_1}\nabla v^\la\nabla\frac{dv^\la}{d\la}-\int_{B_1}|u^\la|^{p-1}u^\la\frac{du^\la}{d\la}\\
&=\int_{B_1}\nabla v^\la\nabla\frac{dv^\la}{d\la}+\int_{B_1}\Delta w^\la\frac{du^\la}{d\la}.
\endaligned\ee
Next, we calculate the term $\int_{B_1}\nabla v^\la\nabla\frac{dv^\la}{d\la}$:

\be\label{10overE1}\aligned
\int_{B_1}\nabla v^\la\nabla\frac{dv^\la}{d\la}&=\int_{\pa B_1}\frac{\pa v^\la}{\pa r}\frac{dv^\la}{d\la}-\int_{B_1}\Delta v^\la\frac{dv^\la}{d\la}\\
&=\int_{\pa B_1}\frac{\pa v^\la}{\pa r}\frac{dv^\la}{d\la}-\int_{B_1}w^\la\Delta \frac{du^\la}{d\la}\\
&=\int_{\pa B_1}\frac{\pa v^\la}{\pa r}\frac{dv^\la}{d\la}-\int_{ \pa B_1}w^\la\frac{\pa}{\pa r}\frac{du^\la}{d\la}+\int_{  B_1}\nabla w^\la\nabla\frac{du^\la}{d\la}\\
&=\int_{\pa B_1}\frac{\pa v^\la}{\pa r}\frac{dv^\la}{d\la}-\int_{\pa B_1}w^\la\frac{\pa}{\pa r}\frac{du^\la}{d\la}\\
&\quad+\int_{\pa B_1}\frac{\pa w^\la}{\pa r}\frac{du^\la}{d\la}-\int_{B_1}\Delta w^\la\frac{du^\la}{d\la}.
\endaligned\ee
By \eqref{10overE} and \eqref{10overE1} we obtain that
\be\label{10overee}\aligned
\frac{d}{d\la}\overline{E}(u^\la,1)=&\int_{\pa B_1}\frac{\pa v^\la}{\pa r}\frac{dv^\la}{d\la}
+\int_{\pa B_1}\frac{\pa w^\la}{\pa r}\frac{du^\la}{d\la}-\int_{\pa B_1}w^\la\frac{\pa}{\pa r}\frac{du^\la}{d\la}.
\endaligned\ee
Recalling  \eqref{10uvw} and
differentiating it with respect to $\la$, we have
\be\nonumber\aligned
\frac{du^\la(x)}{d\la}&=\frac{1}{\la}\big(\frac{6}{p-1}u^\la(x)+r\pa_r u^\la(x)\big),\\
\frac{dv^\la(x)}{d\la}&=\frac{1}{\la}\big((\frac{6}{p-1}+2)v^\la(x)+r\pa_r v^\la(x)\big),\\
\frac{dw^\la(x)}{d\la}&=\frac{1}{\la}\big((\frac{6}{p-1}+4)w^\la(x)+r\pa_r w^\la(x)\big).\\
\endaligned
\ee
Differentiating the above equations with respect to $\la$ again we get
\be\nonumber
\la\frac{d^2u^\la(x)}{d\la^2}+\frac{du^\la(x)}{d\la}=\frac{6}{p-1}\frac{du^\la(x)}{d\la}+r\pa_r\frac{du^\la}{d\la}.
\ee
Hence, for $x\in B_1$, we have
\be\nonumber\aligned
\pa_r(u^\la(x))&=\la\frac{du^\la}{d\la}-\frac{6}{p-1}u,\\
\pa_r(\frac{du^\la(x)}{d\la})&=\la\frac{d^2u_e^\la(x)}{d\la^2}+(1-\frac{6}{p-1})\frac{du^\la}{d\la},\\
\pa_r(v^\la(x))&=\la\frac{dv^\la}{d\la}-(\frac{6}{p-1}+2)v^\la,\\
\pa_r(w^\la(x))&=\la\frac{d w^\la}{d\la}-(\frac{6}{p-1}+4)w^\la.\\
\endaligned
\ee
Plugging these equations into \eqref{10overee}, we get that
\be\label{10overlineE}\aligned
\frac{d}{d\la}\overline{E}(u^\la,1)&=\int_{\pa B_1}\big(\la\frac{dv^\la}{d\la}\frac{dv^\la}{d\la}-(\frac{6}{p-1}+2)v^\la\frac{dv^\la}{d\la}\big)\\
&\quad+\big(\la\frac{d w^\la}{d\la}\frac{u^\la}{d\la}-(\frac{6}{p-1}+4)w^\la\frac{du^\la}{d\la}\big)\\
&\quad -\big(\la w^\la\frac{d^2 u^\la}{d\la^2}+(1-\frac{6}{p-1}w^\la\frac{du^\la}{d\la}\big)\\
&=\underbrace{\int_{ \pa B_1}\big[\la\frac{dv^\la}{d\la}\frac{dv^\la}{d\la}-(\frac{6}{p-1}+2)v^\la\frac{dv^\la}{d\la}\big]}\\
&\quad+\underbrace{\big[\la \frac{d w^\la}{d\la}\frac{du^\la}{d\la}-\la w^\la\frac{d^2u^\la}{d\la^2}\big]-5w^\la\frac{du^\la}{d\la}}.\\
&:=\overline{E}_{d1}(u^\la,1)+\overline{E}_{d2}(u^\la,1)
\endaligned
\ee

\subsection{The calculations  of $\frac{\pa^j}{\pa r^j}u^\la$  and  $\frac{\pa^i}{\pa \la^i}u^\la,\; i,j=1,2,3,4$}

Note
\be\label{10ue0}
\la\frac{du^\la}{d\la}=\frac{6}{p-1}u^\la+r\frac{\pa}{\pa r}u^\la.
\ee
Differentiating \eqref{10ue0} once, twice and  thrice with respect to $\la$ respectively,  we have
\be\label{10ue1}
\la\frac{d^2u^\la}{d\la^2}+\frac{du^\la}{d\la}=\frac{6}{p-1}\frac{du^\la}{d\la}+r\frac{\pa}{\pa r}\frac{du^\la}{d\la},
\ee
\be\label{10ue2}
\la\frac{d^3u^\la}{d\la^3}+2\frac{d^2u^\la}{d\la^2}=\frac{6}{p-1}\frac{d^2u^\la}{d\la^2}+r\frac{\pa}{\pa r}\frac{d^2u^\la}{d\la^2},
\ee
\be\label{10ue3}
\la\frac{d^4u^\la}{d\la^4}+3\frac{d^3u^\la}{d\la^3}=\frac{6}{p-1}\frac{d^3u^\la}{d\la^3}+r\frac{\pa}{\pa r}\frac{d^3u^\la}{d\la^3}.
\ee
Similarly, differentiating \eqref{10ue0} once, twice and  thrice with respect to $r$ respectively we have
\be\label{10ue4}
\la\frac{\pa}{\pa r}\frac{du^\la}{d\la}=(\frac{6}{p-1}+1)\frac{\pa}{\pa r}u^\la+r\frac{\pa^2}{\pa r^2}u^\la,
\ee
\be\label{10ue5}
\la\frac{\pa^2}{\pa r^2}\frac{du^\la}{d\la}=(\frac{6}{p-1}+2)\frac{\pa^2}{\pa r^2}u^\la+r\frac{\pa^3}{\pa r^3}u^\la,
\ee
\be\label{10ue6}
\la\frac{\pa^3}{\pa r^3}\frac{du^\la}{d\la}=(\frac{6}{p-1}+3)\frac{\pa^3}{\pa r^3}u^\la+r\frac{\pa^4}{\pa r^4}u^\la.
\ee
From \eqref{10ue0}, on $\pa B_1$, we have
\be\nonumber
\frac{\pa u^\la}{\pa r}=\la \frac{du^\la}{d\la}-\frac{6}{p-1}u^\la.
\ee
Next from \eqref{10ue1}, on $\pa B_1$,  we derive that
\be\nonumber
\frac{\pa}{\pa r}\frac{d u^\la}{d\la}=\la\frac{d^2u^\la}{d\la^2}+(1-\frac{6}{p-1})\frac{du^\la}{d\la}.
\ee
From \eqref{10ue4}, combining the two equations above, on $\pa B_1$, we get
\be\label{10r31}\aligned
\frac{\pa^2}{\pa r^2} u^\la&=\la\frac{\pa}{\pa r}\frac{du^\la}{d\la}-(1+\frac{6}{p-1})\frac{\pa}{\pa r} u^\la\\
&=\la^2\frac{d^2 u^\la}{d\la^2}-\la\frac{12}{p-1}\frac{du^\la}{d\la}+(1+\frac{6}{p-1})\frac{6}{p-1}u^\la.
\endaligned\ee
Differentiating \eqref{10ue1} with respect to $r$, and combine with \eqref{10ue1} and \eqref{10ue2}, we get that
\be\label{10r32}\aligned
\frac{\pa^2}{\pa r^2}\frac{du^\la}{d\la}&=\la\frac{\pa}{\pa r}\frac{d^2 u^\la}{d\la^2}-\frac{6}{p-1}\frac{\pa}{\pa r}\frac{du^\la}{d\la}\\
&=\la^2\frac{d^3 u^\la}{d\la^3}+(2-\frac{12}{p-1})\la\frac{d^2 u^\la}{d\la^2}-(1-\frac{6}{p-1})\frac{6}{p-1}\frac{du^\la}{d\la}.
\endaligned\ee
From \eqref{10ue5}, on $\pa B_1$, combine with \eqref{10r31} and \eqref{10r32}, we have
\be\label{10ue7}\aligned
\frac{\pa^3}{\pa r^3}u^\la=&\la\frac{\pa^2}{\pa r^2}\frac{du^\la}{d\la}-(2+\frac{6}{p-1})\frac{\pa^2}{\pa r^2}u^\la\\
=&\la^3\frac{d^3 u^\la}{d\la^3}-\la^2\frac{18}{p-1}\frac{d^2 u^\la}{d\la^2}+\la(\frac{18}{p-1}+\frac{108}{(p-1)^2})\frac{du^\la}{d\la}\\
&-(2+\frac{6}{p-1})(1+\frac{6}{p-1})\frac{6}{p-1}u^\la.
\endaligned\ee
Now differentiating \eqref{10ue1} once with respect to $r$, we get
\be\nonumber
\la\frac{\pa^2}{\pa r^2}\frac{d^2 u^\la}{d\la^2}=(\frac{6}{p-1}+1)\frac{\pa^2}{\pa r^2}\frac{du^\la}{d\la}+r\frac{\pa^3}{\pa r^3}\frac{du^\la}{d\la},
\ee
then on $\pa B_1$, we have
\be\label{10ue8}
\frac{\pa^3}{\pa r^3}\frac{du^\la}{d\la}=\la\frac{\pa^2}{\pa r^2}\frac{d^2 u^\la}{d\la^2}-(\frac{6}{p-1}+1)\frac{\pa^2}{\pa r^2}\frac{du^\la}{d\la}.
\ee
Now differentiating \eqref{10ue2} twice with respect to $r$, we get
\be\nonumber
\la\frac{\pa}{\pa r}\frac{d^3 u^\la}{d\la^3}=(\frac{6}{p-1}-1)\frac{\pa}{\pa r}\frac{d^2 u^\la}{d\la^2}+r\frac{\pa^2}{\pa r^2}\frac{d^2 u^\la}{d\la^2}.
\ee
Hence on $\pa B_1$, combining with \eqref{10ue2} and \eqref{10ue3} there holds
\be\label{10ue9}
\aligned
\frac{\pa^2}{\pa r^2}&\frac{d^2 u^\la}{d\la^2}=\la\frac{\pa}{\pa r}\frac{d^3 u^\la}{d\la^3}+(1-\frac{6}{p-1})\frac{\pa}{\pa r}\frac{d^2 u^\la}{d\la^2}\\
=&\la^2\frac{d^4 u^\la}{d\la^4}+\la(4-\frac{12}{p-1})\frac{d^3 u^\la}{d\la^3}+(1-\frac{6}{p-1})(2-\frac{6}{p-1})\frac{d^2 u^\la}{d\la^2}.
\endaligned\ee
Now differentiating \eqref{10ue1} with respect to $r$, we have
\be\nonumber
\la\frac{\pa}{\pa r}\frac{d^2 u^\la}{d\la^2}=\frac{6}{p-1}\frac{\pa}{\pa r}\frac{du^\la}{d\la}+r\frac{\pa^2}{\pa r^2}\frac{du^\la}{d\la}.
\ee
From this combined with \eqref{10ue1} and \eqref{10ue2}, on $\pa B_1$, we have
\be\label{10ue10}\aligned
\frac{\pa^2}{\pa r^2}\frac{du^\la}{d\la}&=\la\frac{\pa}{\pa r}\frac{d^2 u^\la}{d\la^2}-\frac{6}{p-1}\frac{\pa}{\pa r}\frac{d u^\la}{d\la}\\
&=\la^2\frac{d^3 u^\la}{d\la^3}+\la(2-\frac{12}{p-1})\frac{d^2 u^\la}{d\la^2}-\frac{6}{p-1}(1-\frac{6}{p-1})\frac{du^\la}{d\la}.
\endaligned\ee
Now from \eqref{10ue8}, combining with \eqref{10ue9} and \eqref{10ue10}, we get
\be\label{10ue11}\aligned
\frac{\pa^3}{\pa r^3}\frac{du^\la}{d\la}=&\la^3\frac{d^4 u^\la}{d\la^4}+\la^2(3-\frac{18}{p-1})\frac{d^3 u^\la}{d\la^3}-\la(1-\frac{6}{p-1})\frac{18}{p-1}\frac{d^2 u^\la}{d\la^2}\\
&+(1-\frac{6}{p-1})(1+\frac{6}{p-1})\frac{6}{p-1}\frac{du^\la}{d\la}.
\endaligned\ee
From \eqref{10ue6}, on $\pa B_1$, combining with \eqref{10ue11} yields
\be\nonumber\aligned
\frac{\pa^4}{\pa r^4}u^\la=&\la \frac{\pa^3}{\pa r^3}\frac{du^\la}{d\la}-(3+\frac{6}{p-1})\frac{\pa^3}{\pa r^3}u^\la\\
=&\la^4\frac{d^4 u^\la}{d\la^4}-\la^3\frac{24}{p-1}\frac{d^3 u^\la}{d\la^3}+\la^2(2+\frac{12}{p-1})\frac{18}{p-1}\frac{d^2 u^\la}{d\la^2}\\
&-\la (1+\frac{6}{p-1})(1+\frac{3}{p-1})\frac{48}{p-1}\frac{du^\la}{d\la}\\
&+(3+\frac{6}{p-1})(2+\frac{6}{p-1})(1+\frac{6}{p-1})\frac{6}{p-1}u^\la.\\
\endaligned\ee
In summary, we have that
\be\nonumber\aligned
\frac{\pa^3}{\pa r^3}u^\la
=&\la^3\frac{d^3 u^\la}{d\la^3}-\la^2\frac{18}{p-1}\frac{d^2 u^\la}{d\la^2}+\la(\frac{18}{p-1}+\frac{108}{(p-1)^2})\frac{du^\la}{d\la}\\
&-(2+\frac{6}{p-1})(1+\frac{6}{p-1})\frac{6}{p-1}u^\la
\endaligned\ee
and
\be\nonumber
\frac{\pa ^2}{\pa r^2}u^\la=\la^2\frac{d^2 u^\la}{d\la^2}-\la\frac{12}{p-1}\frac{du^\la}{d\la}+(1+\frac{6}{p-1})\frac{6}{p-1}u^\la
\ee
\be\nonumber
\frac{\pa u^\la}{\pa r}=\la \frac{du^\la}{d\la}-\frac{6}{p-1}u^\la.
\ee
\subsection{On the operator $\Delta^2$  and its representation }

Note that
\be\nonumber\aligned
\Delta u=&\nabla\cdot(\nabla u)
=&u_{rr}+\frac{n-1}{r}u_r+\frac{1}{r^2}\mathbf{ div} _{\q}(\nabla_{\q}u).
\endaligned\ee
Set $v:=\Delta u$  and $w:=\Delta^2 u$.
Then
\be\nonumber\aligned
w=&\Delta v=v_{rr}+\frac{n-1}{r}v_r+\frac{1}{r^2} \mathbf{ div} _{\q}(\nabla_{\q}v)\\
=&\pa_{rrrr}u+\frac{2(n-1)}{r}\pa_{rrr}u+\frac{(n-1)(n-3)}{r^2}\pa_{rr}u-\frac{(n-1)(n-3)}{r^3}\pa_r u\\
&+r^{-4}{\bf div}_{\q}(\nabla_\theta(\mathbf{ div}_{\q}(\nabla_{\q}u))\\
&+2r^{-2}{\bf div}_{\q}(\nabla_{\q}(u_{rr}+\frac{n-3}{r}u_r))\\
&-2(n-4)r^{-4}\mathbf{ div}_{\q}(\nabla_{\q}u).\\
\endaligned\ee
On $\pa B_1$, we have
\be\nonumber\aligned
w=&\underbrace{\pa_{rrrr}u+2(n-1)\pa_{rrr}u+(n-1)(n-3)\pa_{rr}u-(n-1)(n-3)\pa_r u}\\
&\underbrace{+\mathbf{ div}_{\q}(\nabla_{\q}(\mathbf{div}_{\q}(\nabla_{\q}u))}\\
&\underbrace{+2\mathbf{div}_{\q}(\nabla_{\q}(u_{rr}+\frac{n-3}{r}u_r))}\\
&\underbrace{-2(n-4)\mathbf{div}_{\q}(\nabla_{\q}u)}\\
&:=I(u)+J(u)+K(u)+L(u).\\
\endaligned\ee
By these notations, we can rewrite the term $\overline{E}_{d2}(u^\la,1)$ appeared  in \eqref{10overlineE} as following
\be\label{10ijk}\aligned
&\overline{E}_{d2}(u^\la,1)\\
&=\int_{\pa B_1} (\la\frac{d w^\la}{d\la}\frac{d u^\la}{d\la}-\la w^\la\frac{d^2 u^\la}{d\la^2})-5 w^\la\frac{d u^\la}{d\la}\\
&=\underbrace{\int_{\pa B_1}\la\frac{d}{d\la}I(u^\la)\frac{d u^\la}{d\la}-\la I(u^\la)\frac{d^2 u^\la}{d\la^2}-5I(u^\la)\frac{du^\la}{d\la}}\\
&\quad \underbrace{+\int_{\pa B_1}\la\frac{d}{d\la}J(u^\la)\frac{d u^\la}{d\la}-\la J(u^\la)\frac{d^2 u^\la}{d\la^2}-5J(u^\la)\frac{du^\la}{d\la}}\\
&\quad \underbrace{+\int_{\pa B_1}\la\frac{d}{d\la}K(u^\la)\frac{d u^\la}{d\la}-\la K(u^\la)\frac{d^2 u^\la}{d\la^2}-5K(u^\la)\frac{du^\la}{d\la}}\\
&\quad \underbrace{+\int_{\pa B_1}\la\frac{d}{d\la}L(u^\la)\frac{du^\la}{d\la}-\la L(u^\la)\frac{d^2u^\la}{d\la^2}-5L(u^\la)\frac{du^\la}{d\la}},
\endaligned\ee
correspondingly, we  rewrite $\overline{E}_{d2}(u^\la,1)$ as
\be\nonumber\aligned
\overline{E}_{d2}(u^\la,1):&=\mathcal{I}+\mathcal{J}+\mathcal{K}+\mathcal{L}\\
:&=I_1+I_2+I_3+J_1+J_2+J_3+K_1+K_2+K_3+L_1+L_2+L_3,
\endaligned\ee
where $I_1,I_2,I_3,J_1,J_2,J_3,K_1,K_2,K_3,L_1,L_2,L_3 $ successively corresponding   to the $12$ terms in \eqref{10ijk}.
By the conclusions of subsection $2.2$, we have
\be\label{10iuiu}\aligned
I(u^\la)&=\pa_{rrrr} u^\la+2(n-1)\pa_{rrr} u^\la\\
&\quad +(n-1)(n-3)\pa_{rr} u^\la-(n-1)(n-3)\pa_r u^\la\\
&=\la^4\frac{d^4 u^\la}{d\la^4}+\la^3\big(2(n-1)-\frac{24}{p-1})\big)\frac{d^3 u^\la}{d\la^3}\\
&\quad+\la^2\big[\frac{36}{p-1}(1+\frac{6}{p-1})-(n-1)\frac{36}{p-1}+(n-1)(n-3)\big]\frac{d^2 u^\la}{d\la^2}\\
&\quad+\la\big[-\frac{24}{p-1}(1+\frac{6}{p-1})(2+\frac{6}{p-1})+2(n-1)\frac{18}{p-1}(1+\frac{6}{p-1})\\
&\quad+(n-1)(n-3)(-\frac{12}{p-1}-1)\big]\frac{d u^\la}{d\la}\\
&\quad+\big[(1+\frac{6}{p-1})(2+\frac{6}{p-1})(3+\frac{6}{p-1})\frac{6}{p-1}\\
&\quad-(n-1)(1+\frac{6}{p-1})(2+\frac{6}{p-1})\frac{12}{p-1}\\
&\quad+(n-1)(n-3)(\frac{6}{p-1}+2)\frac{6}{p-1}\big]u^\la.\\
\endaligned\ee
For convenience, we denote that
\be\label{10iue}
I(u^\la)=\la^4\frac{d^4 u^\la}{d\la^4}+\la^3 \delta_1\frac{d^3 u^\la}{d\la^3}+\la^2\delta_2\frac{d^2 u^\la}{d\la^2}+\la\delta_3\frac{du^\la}{d\la}+\delta_4 u^\la,
\ee
where $\delta_i$ are the corresponding coefficients of $\la^i\frac{d^i u^\la}{d\la^i}$ appeared in \eqref{10iuiu} for $i=1,2,3,4$.
Now taking the derivative of \eqref{10iue} with respect to $\la$ we get
\be\label{10diue}\aligned
\frac{d}{d\la}I(u^\la)=&\la^4\frac{d^5 u^\la}{d\la^5}+\la^3 (\delta_1+4)\frac{d^4 u^\la}{d\la^4}+\la^2(3\delta_1+\delta_2)\frac{d^3 u^\la}{d\la^3}\\
&+\la(2\delta_2+\delta_3)\frac{d^2u^\la}{d\la^2}+(\delta_3+\delta_4)\frac{du^\la}{d\la}.
\endaligned\ee
Since
\be\label{10urr}\aligned
\pa_{rr}& u^\la+(n-3)\pa_r u^\la \\ =&\la^2\frac{d^2u^\la}{d\la^2}+\la(n-3-\frac{12}{p-1})\frac{du^\la}{d\la}+\frac{6}{p-1}(4+\frac{6}{p-1}-n)u^\la\\
:=&\la^2\frac{d^2 u^\la}{d\la^2}+\la\alpha\frac{d^\la}{d\la}+\beta u^\la.
\endaligned\ee
Hence,
\be\label{8urrla}\aligned
\frac{d}{d\la}[\pa_{rr}& u^\la+(n-1)\pa_r u^\la]=\la^2\frac{d^3 u^\la}{d\la^3}+\la(\alpha+2)\frac{d^2 u^\la}{d\la^2}+(\alpha+\beta)\frac{du^\la}{d\la},
\endaligned\ee
here $\alpha=n-3-\frac{12}{p-1}$ and $\beta=\frac{6}{p-1}(4+\frac{6}{p-1}-n)$.

\subsection{The computations  of $I_1,I_2,I_3$ and $\mathcal{I}$}
We start with
\be\label{10i1}\aligned
I_1:=&\int_{\pa B_1}\la \frac{d}{d\la}I(u^\la)\frac{du^\la}{d\la}\\
=&\int_{\pa B_1}\big(\la^5\frac{d^5 u^\la}{d\la^5}+\la^4(4+\delta_1)\frac{d^4 u^\la}{d\la^4}+\la^3(3\delta_1+\delta_2)\frac{d^3 u^\la}{d\la^3}\\
&+\la^2(2\delta_2+\delta_3)\frac{d^2 u^\la}{d\la^2}+\la(\delta_3+\delta_4)\frac{du^\la}{d\la}\big)\frac{du^\la}{d\la}\\
=&\frac{d}{d\la}\int_{\pa B_1}\Big[\la^5\frac{d^4 u^\la}{d\la^4}\frac{d u^\la}{d\la}-\la^5\frac{d^3 u^\la}{d\la^3}\frac{d^2 u^\la}{d\la^2}+(\delta_1-1)\la^4\frac{d^3 u^\la}{d\la^3}\frac{d u^\la}{d\la}\\
&+(4-\delta_1+\delta_2)\la^3\frac{d^2 u^\la}{d\la^2}\frac{d u^\la}{d\la}+\frac{3\delta_1-\delta_2+\delta_3-12}{2}\la^2(\frac{d u^\la}{d\la})^2\Big]\\
&+\int_{\pa B_1} \Big[(12-3\delta_1+\delta_2+\delta_4)\la(\frac{d u^\la}{d\la})^2\\
&+(\delta_1-4-\delta_2)\la^3(\frac{d^2 u^\la}{d\la^2})^2+\la^5(\frac{d^3 u^\la}{d\la^3})^2\Big]\\
&+\int_{\pa B_1}(6-\delta_1)\la^4\frac{d^3 u^\la}{d\la^3}\frac{d^2 u^\la}{d\la^2},
\endaligned
\ee
where $\delta_i,i=1,2,3,4, $ are defined in \eqref{10iuiu} and \eqref{10iue}. In this computation, we denote that $f=u^\la,f':=\frac{du^\la}{d\la}$ and  we have used the fact that
\be\nonumber\aligned
\la^5 f'''''f'=&\big[\la^5f''''f'-\la^5f'''f''-5\la^4f'''f'+20\la^3f''f'-30\la^2f'f'\big]'\\
&+60\la(f')^2-20\la^3(f'')^2+\la^5(f''')^2+10\la^4f'''f'',
\endaligned\ee
\be\nonumber\aligned
\la^4f''''f'=\big[\la^4f'''f'-4\la^3f''f'+6\la^2f'f'\big]'-12\la(f')^2+4\la^3(f'')^2-\la^4f'''f'',
\endaligned\ee
\be\nonumber\aligned
\la^3f'''f'=\big[\la^3f''f'-\frac{3\la^2}{2}f'f'\big]'+3\la(f')^2-\la^3(f'')^2,
\endaligned\ee
and
\be\nonumber\aligned
\la^2 f''f'=\big[\frac{\la^2}{2}f'f'\big]'-\la(f')^2.
\endaligned\ee

\be\label{8i2}\aligned
&I_2:=-\la\int_{\pa B_1} I(u_e^\la)\frac{d^2 u^\la}{d\la^2}\\
&=-\la\int_{\pa B_1} \big(\la^4\frac{d^4 u^\la}{d\la^4}+\la^3 \delta_1\frac{d^3 u_e^\la}{d\la^3}+\la^2\delta_2\frac{d^2 u^\la}{d\la^2}+\la\delta_3\frac{du^\la}{d\la}+\delta_4 u^\la\big)\frac{d^2 u^\la}{d\la^2}\\
&=\frac{d}{d\la}\int_{\pa B_1}\big[-\la^5\frac{d^3 u^\la}{d\la^3}\frac{d^2 u^\la}{d\la^2}-\delta_4\la\frac{d u^\la}{d\la}u^\la\big]\\
&\quad+\int_{\pa B_1}\big[\la^5(\frac{d^3 u^\la}{d\la^3})^2-\delta_2\la^3(\frac{d^2 u^\la}{d\la^2})^2+\delta_4\la(\frac{d u^\la}{d\la})^2\big]\\
&\quad +\int_{\pa B_1}\big[(5-\delta_1)\la^4\frac{d^3 u^\la}{d\la^3}\frac{d^2 u^\la}{d\la^2}-\delta_3\la^2\frac{d^2 u^\la}{d\la^2}\frac{d u^\la}{d\la}+\delta_4\frac{d u^\la}{d\la}u^\la\big],
\endaligned\ee
here we have used that
\be\nonumber\aligned
-\la^5f''''f''=\big[-\la^5f'''f''\big]'+5\la^4f'''f''+\la^5(f''')^2
\endaligned\ee
and
\be\nonumber\aligned
-\la f''f=[-\la f'f]'+f'f+\la(f')^2.
\endaligned\ee
Further,
\be\label{8iu3}\aligned
I_3:=&-5\int_{\pa B_1} I(u^\la)\frac{d u^\la}{d\la}\\
=&-5\int_{\pa B_1}\big[\la^4\frac{d^4 u^\la}{d\la^4}+\la^3 \delta_1\frac{d^3 u^\la}{d\la^3}+\la^2\delta_2\frac{d^2 u^\la}{d\la^2}+\la\delta_3\frac{du^\la}{d\la}+\delta_4 u^\la\big]\frac{d u^\la}{d\la}\\
=&\frac{d}{d\la}\int_{\pa B_1}\big[-5\frac{d^3 u^\la}{d\la^3}\frac{d u^\la}{d\la}+(20-5\delta_1)\la^3\frac{d^2 u^\la}{d\la^2}\frac{d u^\la}{d\la}\big]\\
&+\int_{\pa B_1}\big[(5\delta_1-20)\la^3(\frac{d^2 u^\la}{d\la^2})^2-5\delta_3\la(\frac{d u^\la}{d\la})^2\big]\\
&+\int_{\pa B_1}\big[5\la^4\frac{d^3 u^\la}{d\la^3}\frac{d^2 u^\la}{d\la^2}+(15\delta_1-60-5\delta_2)\frac{d^2 u^\la}{d\la^2}\frac{d u^\la}{d\la}-5\delta_4\frac{d u^\la}{d\la}u^\la\big],
\endaligned\ee
here we have use that
\be\nonumber\aligned
-\la^4 f''''f'=\big[-5\la^4f'''f'+20\la^3f''f'\big]'-20\la^3(f'')^2-60\la^2f''f'+5\la^4f'''f''
\endaligned\ee
and
\be\nonumber\aligned
-\la^3 f'''f'=\big[-\la^3f''f'\big]'+3\la^2 f''f'+\la^3 (f'')^2.
\endaligned\ee
Summing  up $I_1,I_2, I_3$, we can get the term $\mathcal{I}$.
\be\label{8I}\aligned
\mathcal{I}:=& I_1+I_2+I_3\\
=&\frac{d}{d\la}\int_{\pa B_1}\big[\la^5\frac{d^4 u_e^\la}{d\la^4}\frac{d u_e^\la}{d\la}-2\la^5\frac{d^3 u_e^\la}{d\la^3}\frac{d^2 u_e^\la}{d\la^2}\\
&+(\delta_1-6)\la^4\frac{d^3 u^\la}{d\la^3}\frac{d u^\la}{d\la}+(24-6\delta_1+\delta_2)\la^3\frac{d^2 u^\la}{d\la^2}\frac{d u^\la}{d\la}\\
&+\frac{3\delta_1-\delta_2+\delta_3-12}{2}\la^2\frac{d u^\la}{d\la}\frac{d u^\la}{d\la}\\
&+(\delta_1-8)\la^4(\frac{d^2 u^\la}{d\la^2})^2-\delta_4\la\frac{d u^\la}{d\la}u^\la-2\delta_4(u^\la)^2\big]\\
&+\int_{\pa B_1}\big[(10\delta_1-2\delta_2-56)\la^3(\frac{d^2 u^\la}{d\la^2})^2\\
&+(12-3\delta_1+\delta_2-5\delta_3+2\delta_4)\la(\frac{d u^\la}{d\la})^2\big]\\
&+\int_{\pa B_1}\big[(15\delta_1-5\delta_2-\delta_3-60)\la^2\frac{d^2 u^\la}{d\la^2}\frac{d u^\la}{d\la}\big]\\
&+2\la^5\int_{\pa B_1}(\frac{d^3 u^\la}{d\la^3})^2.
\endaligned\ee
Since $u^\la(x)=\la^{\frac{6}{p-1}}u(\la x)$, we have the following
\be\nonumber\aligned
\la^4&\frac{d^4 u^\la}{d\la^4}=\la^{\frac{6}{p-1}}\big[\frac{6}{p-1}(\frac{6}{p-1}-1)(\frac{6}{p-1}-2)(\frac{6}{p-1}-3)u(\la x)\\
&+\frac{24}{p-1}(\frac{6}{p-1}-1)(\frac{6}{p-1}-2)r\la\pa_r u(\la x)\\
&+\frac{36}{p-1}(\frac{6}{p-1}-1)r^2\la^2\pa_{rr}u(\la x)\\
&+\frac{24}{p-1}r^3\la^3\pa_{rrr}u(\la x)+r^4\la^4\pa_{rrrr}u(\la x)\big],\\
\endaligned\ee
and
\be\nonumber\aligned
\la^3\frac{d^3 u^\la}{d\la^3}&=\la^{\frac{6}{p-1}}\big[\frac{6}{p-1}(\frac{6}{p-1}-1)(\frac{6}{p-1}-2)u(\la x)\\
&+\frac{18}{p-1}(\frac{6}{p-1}-1)r\la\pa_r u(\la x)\\
&+\frac{18}{p-1}r^2\la^2 \pa_{rr} u(\la x)+r^3\la^3\pa_{rrr} u(\la x)\big],
\endaligned\ee
\be\nonumber\aligned
&\la^2\frac{d^2 u^\la}{d\la^2}\\
&=\la^{\frac{6}{p-1}}\big[\frac{6}{p-1}(\frac{6}{p-1}-1) u(\la x)+\frac{12}{p-1}r\la\pa_r u(\la x)+r^2\la^2\pa_{rr}u(\la x)\big]
\endaligned\ee
and
\be\nonumber\aligned
\la\frac{d u^\la}{d\la}=\la^{\frac{6}{p-1}}\big[\frac{6}{p-1}u(\la x)+r\la\pa_r u(\la x)\big].
\endaligned\ee
Hence, by scaling we have
\be\nonumber\aligned
\frac{d}{d\la}&\int_{\pa B_1}\la^5\frac{d^4 u^\la}{d\la^4}\frac{d u^\la}{d\la}\\
&=\frac{d}{d\la}\int_{\pa B_\la}\la^{6\frac{p+1}{p-1}-n-5}y^b\big[\frac{6}{p-1}(\frac{6}{p-1}-1)(\frac{6}{p-1}-2)(\frac{6}{p-1}-3)u\\
&+\frac{24}{p-1}(\frac{6}{p-1}-1)(\frac{6}{p-1}-2)\la\pa_r u+\frac{36}{p-1}(\frac{6}{p-1}-1)\la^2\pa_{rr}u\\
&+\frac{24}{p-1}\la^3\pa_{rrr}u+\la^4\pa_{rrrr}u\big]\big[\frac{6}{p-1}u+r\la\pa_r u\big],
\endaligned\ee
further,
\be\nonumber\aligned
\frac{d}{d\la}&\int_{\pa B_1}\la^5\frac{d^3 u^\la}{d\la^3}\frac{d^2 u^\la}{d\la^2}\\
&=\frac{d}{d\la}\int_{\pa B_\la}\la^{6\frac{p+1}{p-1}-n-5}\big[\frac{6}{p-1}(\frac{6}{p-1}-1)(\frac{6}{p-1}-2)u\\
&\quad+\frac{18}{p-1}(\frac{6}{p-1}-1)\la\pa_r u
\quad\quad+\frac{18}{p-1}\la^2 \pa_{rr} u+\la^3\pa_{rrr} u\big]\\
&\quad \quad \big[\frac{6}{p-1}(\frac{6}{p-1}-1) u+\frac{12}{p-1}\la\pa_r u+\la^2\pa_{rr}u\big],
\endaligned\ee

\be\nonumber\aligned
\frac{d}{d\la}&\int_{\pa B_1}\la^4\frac{d^3 u^\la}{d\la^3}\frac{d u^\la}{d\la}\\
&=\frac{d}{d\la}\int_{\pa B_\la}\la^{6\frac{p+1}{p-1}-n-5}\big[\frac{6}{p-1}(\frac{6}{p-1}-1)(\frac{6}{p-1}-2)u\\
&+\frac{18}{p-1}(\frac{6}{p-1}-1)\la\pa_r u
+\frac{18}{p-1}\la^2 \pa_{rr} u+\la^3\pa_{rrr} u\big]\\
&\big[\frac{6}{p-1}u+\la\pa_r u\big].
\endaligned\ee
On the other hand,
\be\nonumber\aligned
\frac{d}{d\la}&\int_{\pa B_1}\la^3\frac{d^2 u^\la}{d\la^2}\frac{d u^\la}{d\la}\\
&=\frac{d}{d\la}\int_{\pa B_\la}\la^{6\frac{p+1}{p-1}-n-5}
\big[\frac{6}{p-1}(\frac{6}{p-1}-1) u\\
&\quad\quad +\frac{12}{p-1}\la\pa_r u
+\la^2\pa_{rr}u][\frac{6}{p-1}u+\la\pa_r u\big],
\endaligned\ee

\be\nonumber\aligned
\frac{d}{d\la}\int_{\pa B_1}\la^2\frac{d u^\la}{d\la}\frac{d u^\la}{d\la}
=\frac{d}{d\la}\int_{\pa B_\la}\la^{6\frac{p+1}{p-1}-n-5}
\big[\frac{6}{p-1}u+\la\pa_r u\big]^2,
\endaligned\ee

\be\nonumber\aligned
\frac{d}{d\la}\int_{\pa B_1}\la^4\frac{d^2 u^\la}{d\la^2}\frac{d^2 u^\la}{d\la^2}
&=\frac{d}{d\la}\int_{\pa B_\la}\la^{6\frac{p+1}{p-1}-n-5}
\Big(\frac{6}{p-1}(\frac{6}{p-1}-1) u\\
&\quad\quad+\frac{12}{p-1}\la\pa_r u+\la^2\pa_{rr}u\Big)^2,
\endaligned\ee

\be\nonumber\aligned
\frac{d}{d\la}\int_{\pa B_1}\la\frac{d u^\la}{d\la}u^\la
=\frac{d}{d\la}\int_{\pa B_\la}\la^{6\frac{p+1}{p-1}-n-5}
\big[\frac{6}{p-1}u)+\la\pa_r u)\big]u,
\endaligned\ee
and
\be\nonumber\aligned
\frac{d}{d\la}&\int_{\pa B_1} u^\la=\frac{d}{d\la}\int_{\pa B_\la}\la^{6\frac{p+1}{p-1}-n-5}
u^2.
\endaligned\ee

\subsection{The computations of $J_i,K_i,L_i (i=1,2,3)$ and $\mathcal{J},\mathcal{K},\mathcal{L}$}
We begin with
\be\label{10j1}\aligned
J_1:=&\int_{\pa B_1}\la \frac{d}{d\la}J(u^\la)\frac{du^\la}{d\la}=\int_{\pa B_1}\la J(\frac{du^\la}{d\la})\frac{du^\la}{d\la}\\
=&\la\int_{\pa B_1} \mathbf{div}_{\q}\big(\nabla_{\q}(\mathbf{div}_{\q}(\nabla_{\q}\frac{du^\la}{d\la}))\big)\frac{du^\la}{d\la}\\
=&-\la\int_{\pa B_1}\nabla_{\q}\big( \mathbf{div}_{\q}( \nabla_{\q}\frac{du^\la}{d\la})\big)\nabla_{\q}\frac{du^\la}{d\la}\\
=&\la\int_{\pa B_1}\big[\mathbf{div}_{\q}(\nabla_{\q}\frac{du^\la}{d\la})\big]^2\\
=&\la\int_{\pa B_1}\big[\frac{d}{d\la}\mathbf{div}_{\q}(\nabla_{\q}u^\la)\big]^2.\\
 \endaligned\ee
Here we have used integrating by part formula on the unit  sphere $S^n$. Next
\be\label{10j2}\aligned
J_2:=&-\la\int_{\pa B_1} J(u^\la)\frac{d^2 u^\la}{d\la^2}\\
=&-\la\int_{\pa B_1}\mathbf{div}_{\q}\big(\nabla_{\q}(\mathbf{div}_{\q}(\nabla_{\q}u^\la))\big)\frac{d^2 u^\la}{d\la^2}\\
=&\la\int_{\pa B_1}\nabla_{\theta}(\mathbf{ div}_{\q}(\nabla_{\q}u^\la)\nabla_{\q}\frac{d^2u^\la}{d\la^2}\\
=&-\la\int_{\pa B_1} \mathbf{div}_{\q}( \nabla_{\q}u^\la) \frac{d^2}{d\la^2}\mathbf{div}_{\q}(\nabla_{\q}u^\la)\\
=&\frac{d}{d\la}\int_{\pa B_1}-\la \big[\mathbf{div}_{\q}(\nabla_{\q}u^\la)\big]\frac{d}{d\la}\big[\mathbf{div}_{\q}(\nabla_{\q}u^\la)\big]\\
&\quad+\int_{\pa B_1}\mathbf{ div}_{\q}(\nabla_{\q}u^\la)\cdot \frac{d}{d\la}\mathbf{div}_{\q}(\nabla_{\q}u^\la)\\
&\quad+\la\int_{\pa B_1}\big[\frac{d}{d\la}\mathbf{div}_{\q}( \nabla_{\q}u^\la)\big]^2.
 \endaligned\ee
Here we denote that $g= \mathbf{div}_{\q}(\nabla_{\q}u^\la),g'=\frac{d}{d\la}\mathbf{div}_{\q}(\nabla_{\q}u^\la)$
and we have used the fact that
\be\nonumber
-\la g g''=\big[-\la g g'\big]'+g g'+\la (g')^2=\big[-\la g g'+\frac{1}{2}g^2\big]'+\la(g')^2.
\ee
Furthermore,
\be\label{10j3}\aligned
J_3:=&-5\int_{\pa B_1} J(u^\la)\frac{du^\la}{d\la}\\
=&-5\int_{\pa B_1}\mathbf{div}_{\q}\big(\nabla_{\q}(\mathbf{div}_{\q}(\nabla_{\q}u))\big)\frac{du^\la}{d\la}\\
=&5\int_{\pa B_1} \nabla_{\q}\big(\mathbf{div}_{\q}(\nabla_{\q}u)\big)\nabla_{\q}\frac{d u^\la}{d\la}\\
=&-5\int_{\pa B_1}\mathbf{ div}_{\q}(\nabla_{\q}u)\frac{d}{d\la}\mathbf{div}_{\q}(\nabla_{\q}u).\\
 \endaligned\ee
Therefore, combining with \eqref{10j1}, \eqref{10j2} and \eqref{10j3}, we get that
\be\label{8j}\aligned
\mathcal{J}:=&J_1+J_2+J_3\\
=&2\la\int_{\pa B_1}\big[\frac{d}{d\la}\mathbf{div}_{\q}(\nabla_{\q}u^\la)\big]^2\\
&\quad-4\int_{\pa B_1}\mathbf{div}_{\q}(\nabla_{\q}u)\frac{d}{d\la}\mathbf{div}_{\q}(\nabla_{\q}u)\\
&\quad+\frac{d}{d\la}\int_{\pa B_1}-\la\big[\mathbf{div}_{\q}(\nabla_{\q}u^\la)\big]\frac{d}{d\la}\big[\mathbf{div}_{\q}( \nabla_{\q}u^\la)\big]\\
&=2\la\int_{\pa B_1}\big[\frac{d}{d\la}\mathbf{div}_{\q}(\nabla_{\q}u^\la)\big]^2\\
&\quad-2\frac{d}{d\la}\int_{\pa B_1}\mathbf{div}_{\q}(\nabla_{\q}u^\la)\mathbf{div}_{\q}(\nabla_{\q}u^\la)\\
&\quad+\frac{d}{d\la}\int_{\pa B_1}-\la\big[\mathbf{div}_{\q}(\nabla_{\q}u^\la)\big]\frac{d}{d\la}\big[\mathbf{div}_{\q}( \nabla_{\q}u^\la)\big].
\endaligned\ee
Hence, we get that
\be\label{10jj}\aligned
\mathcal{J}\geq&-2\frac{d}{d\la}\int_{\pa B_1}\mathbf{div}_{\q}(\nabla_{\q}u)\mathbf{div}_{\q}(\nabla_{\q}u)\\
&+\frac{d}{d\la}\int_{\pa B_1}-\la \big[\mathbf{div}_{\q}( \nabla_{\q}u^\la)\big]\frac{d}{d\la}\big[\mathbf{div}_{\q}(\nabla_{\q}u) \\
=&-2\frac{d}{d\la}\int_{\pa B_1} [\mathbf{div}_{\q}(\nabla_{\q}u)\big]^2\\
&+\frac{d}{d\la}\int_{\pa B_1}-\la \frac{d}{d\la}\big[\mathbf{div}_{\q}(\nabla_{\q}u)\big]^2.\\
\endaligned\ee
Note that
\be\nonumber\aligned
&\frac{d}{d\la}\int_{\pa B_1}\big[\mathbf{div}_{\q}(\nabla_{\q}u^\la)\big]^2\\
&=\frac{d}{d\la}\int_{\pa B_\la}\la^{6\frac{p+1}{p-1}-n-5}\big(\la^2\Delta u-\la^2\pa_{rr}u-(n-1)\la\pa_r u\big)^2,
\endaligned\ee
and
\be\nonumber\aligned
&\frac{d}{d\la}\int_{\pa B_1}\la \frac{d}{d\la}\big[\mathbf{div}_{\q}(\nabla_{\q}u^\la)\big]^2\\
&=\frac{d}{d\la}\int_{\pa B_\la}\la^{6\frac{p+1}{p-1}-n-4}\frac{d}{d\la}\big(\la^2\Delta u-\la^2\pa_{rr}u-(n-1)\la\pa_r u\big)^2.
\endaligned\ee

Next we compute $K_1,K_2,K_3$ and $\mathcal{K}$.
\be\label{10k1}\aligned
K_1:=&\la\int_{\pa B_1}\frac{d}{d\la}K(u^\la)\frac{du^\la}{d\la}\\
=&2\la\int_{\pa B_1}\mathbf{div}_{\q}\big(\nabla_{\q}(\frac{d}{d\la}(u_{rr}+(n-1)u_r)\big)\frac{du^\la}{d\la}\\
=&2\la\int_{\pa B_1}\mathbf{div}_{\q}\big(\nabla_{\q}(\la^3\frac{d^3 u^\la}{d\la^3}+\la^2(\alpha+2)\frac{d^2 u^\la}{d\la^2}+\la(\alpha+\beta)\frac{d u_e^\la}{d\la})\big)\frac{du^\la}{d\la}\\
=&-2\la\int_{\pa B_1}\nabla_{\q}\big(\la^3\frac{d^3 u^\la}{d\la^3}+\la^2(\alpha+2)\frac{d^2 u^\la}{d\la^2}+\la(\alpha+\beta)\frac{d u^\la}{d\la}\big)\nabla_{\q}\frac{du^\la}{d\la}\\
=&\frac{d}{d\la}\Big(\int_{\pa B_1}-\la^3\frac{d}{d\la}\big(\frac{d}{d\la}\nabla_{\q} u^\la\big)^2\Big)
+(2-2\alpha)\la^2\int_{\pa B_1}(\frac{d}{d\la}\nabla_{\q} u^\la)\\
\cdot &(\frac{d^2}{d\la^2}\nabla_{\q} u^\la)
+2\la^3\int_{\pa B_1}(\frac{d^2}{d\la^2}\nabla_{\q} u^\la)^2\\
&-(2\alpha+2\beta)\la\int_{\pa B_1}(\frac{d}{d\la}\nabla_{\q} u^\la)^2.
\endaligned\ee

Here we denote that $h=\nabla_{\q} u^\la,h'=\frac{d}{d\la}\nabla_{\q} u^\la$, and used that
\be\nonumber
-\la^3h'h'''=\big[-\frac{\la^3}{2}\frac{d}{d\la}(h')^2\big]'+3\la^2 h'h''+\la^3(h'')^2.
\ee

Next,
\be\label{10k2}\aligned
K_2:=&-\la\int_{\pa B_1} K(u^\la)\frac{d^2 u^\la}{d\la^2}\\
=&-2\la\int_{\pa B_1}\mathbf{div}_{\q}(\nabla_{\q}\big(\la^2\frac{d^2 u^\la}{d\la^2}+\la\alpha\frac{d u^\la}{d\la}+\beta u^\la)\big)\frac{d^2 u^\la}{d\la^2}\\
=&2\la\int_{\pa B_1}\nabla_{\q}\big(\la^2\frac{d^2 u^\la}{d\la^2}+\la\alpha\frac{d u^\la}{d\la}+\beta u^\la\big)\nabla_{\q}\frac{d^2 u^\la}{d\la^2}\\
=&\frac{d}{d\la}\int_{\pa B_1}\big[2\beta\la\nabla_{\q}u^\la \frac{d}{d\la}\nabla_{\q}u^\la-\beta(\nabla_{\q}u^\la)^2\big]\\
&+2\la^3\int_{\pa B_1}(\frac{d^2}{d\la^2}\nabla_{\q}u^\la)^2-2\la\beta\int_{\pa B_1}(\frac{d}{d\la}\nabla_{\q}u^\la)^2\\
&+2\la^2\alpha\int_{\pa B_1}\frac{d}{d\la}\nabla_{\q}u^\la\frac{d^2}{d\la^2}\nabla_{\q}u^\la.
\endaligned\ee

Here we have used that
\be\nonumber
2\la h h''=[2\la h h'-h^2]'-2\la(h')^2.
\ee
Further,
\be\label{10k3}\aligned
K_3:=&-5\int_{\pa B_1} K(u^\la)\frac{du^\la}{d\la}\\
=&-10\int_{\pa B_1}\mathbf{div}_{\q}\big(\nabla_{\q}(\la^2\frac{d^2 u^\la}{d\la^2}+\la\alpha\frac{du^\la}{d\la}+\beta u^\la)\big)\frac{du^\la}{d\la}\\
=&10\int_{\pa B_1}\nabla_{\theta}\big(\la^2\frac{d^2 u^\la}{d\la^2}+\la\alpha\frac{du^\la}{d\la}+\beta u^\la\big)\nabla_{\q}\frac{du^\la}{d\la}\\
=&\frac{d}{d\la}\big[5\beta\int_{\pa B_1} \nabla_{\q} u^\la \nabla_{\q} u^\la\big]+10\la\alpha\int_{\pa B_1}\big(\frac{d}{d\la}\nabla_{\q} u^\la\big)^2\\
&+10\la^2\int_{\pa B_1} \frac{d}{d\la}\nabla_{\q} u^\la \frac{d^2}{d\la^2}\nabla_{\q} u^\la.
\endaligned\ee
Now combine with \eqref{10k1}, \eqref{10k2} and \eqref{10k3}, we get that
\be\label{10k}\aligned
\mathcal{K}:=&K_1+K_2+K_3\\
=&\frac{d}{d\la}\int_{\pa B_1}\big[-\la^3\frac{d}{d\la}(\frac{d}{d\la}\nabla_{\q} u^\la)^2\\
& +2\beta\la\nabla_{\q} u^\la\frac{d}{d\la}\nabla_{\q} u^\la+4\beta(\nabla_{\q} u^\la)^2\big]\\
&+4\la^3\int_{\pa B_1}(\frac{d^2}{d\la^2}\nabla_{\q} u^\la)^2+(8\alpha-4\beta)\la\int_{\pa B_1}(\frac{d}{d\la}\nabla_{\q} u^\la)^2\\
&+12\la^2\int_{\pa B_1}\frac{d}{d\la}\nabla_{\q} u^\la\frac{d^2}{d\la^2}\nabla_{\q} u^\la.\\
&=\frac{d}{d\la}\int_{\pa B_1}\big[-\la^3\frac{d}{d\la}(\frac{d}{d\la}\nabla_{\q} u^\la)^2+6\la^2(\nabla_\theta\frac{du^\la}{d\la})^2\\
& +2\beta\la\nabla_{\q} u^\la\frac{d}{d\la}\nabla_{\q} u^\la+4\beta(\nabla_{\q} u^\la)^2\big]\\
&+\int_{\pa B_1}4\la^3(\frac{d^2}{d\la^2}\nabla_{\q} u^\la)^2+(8\alpha-4\beta-12)(\frac{d}{d\la}\nabla_{\q} u^\la)^2
\endaligned\ee

Notice that by scaling we have
\be\nonumber\aligned
\frac{d}{d\la}\int_{\pa B_1}(\nabla_{\q} u^\la)^2
=\frac{d}{d\la}\int_{\pa B_\la}\la^{6\frac{p+1}{p-1}-n-5}\big[\la^2|\nabla u|^2-\la^2|\pa_r u|^2\big],
\endaligned\ee
\be\nonumber\aligned
\frac{d}{d\la}\int_{\pa B_1}\la\frac{d}{d\la}(\nabla_{\q} u^\la)^2=\frac{d}{d\la}\int_{\pa B_\la}\la^{6\frac{p+1}{p-1}-n-4}\frac{d}{d\la}\big[\la^2|\nabla u|^2-\la^2|\pa_r u|^2\big]
\endaligned\ee
and
\be\nonumber\aligned
 \frac{d}{d\la}\int_{\pa B_1}\la^3\frac{d}{d\la}(\frac{d}{d\la}\nabla_{\q} u^\la)^2=&\frac{d}{d\la}\int_{\pa B_\la}\la^{6\frac{p+1}{p-1}-n-2}\frac{d}{d\la}\big[
 \nabla(\frac{6}{p-1}u+\la \pa_r u)^2\\
 &-|\pa_r(\frac{6}{p-1}u+\la \pa_r u)|^2
 \big]^2.
\endaligned\ee
Finally, we compute $\mathcal{L}$.
\be\nonumber\aligned
L_1:&=\int_{\pa B_1}\la\frac{d}{d\la}L(u^\la)\frac{du^\la}{d\la}=-2(n-4)\la\int_{\pa B_1}\mathbf{div}_{\q}(\nabla_{\q}\frac{du^\la}{d\la})\frac{du^\la}{d\la}\\
&=2(n-4)\la\int_{\pa B_1}(\nabla_{\q}\frac{du^\la}{d\la})^2;
\endaligned\ee

\be\nonumber\aligned
L_2:&=\int_{\pa B_1}-\la L(u^\la)\frac{d^2u^\la}{d\la^2}=2(n-4)\la\int_{\pa B_1}\mathbf{div}_{\q}(\nabla_{\q}u^\la)\frac{d^2 u^\la}{d\la^2}\\
&=-2(n-4)\int_{\pa B_1}\la\nabla_{\q}u^\la\frac{d^2}{d\la^2}\nabla_{\q}u^\la\\
&=(4-n)\int_{\pa B_1}\frac{d}{d\la}\big[2\la\nabla_{\q}u^\la\nabla_{\q}\frac{du^\la}{d\la}-(\nabla_{\q}u^\la)^2\big]+2(n-4)\la\int_{\pa B_1}\big[\frac{d}{d\la}\nabla_{\q}u^\la\big]^2;\\
\endaligned\ee

\be\nonumber\aligned
L_3:&=\int_{\pa B_1}-5L(u^\la)\frac{du^\la}{d\la}=10(n-4)\int_{\pa B_1}\mathbf{div}_{\q}(\nabla_{\q}u^\la)\frac{du^\la}{d\la}\\
&=-10(n-4)\int_{\pa B_1}\nabla_{\q}u^\la\nabla_{\q}\frac{du^\la}{d\la}=-5(n-4)\frac{d}{d\la}\int_{\pa B_1}[\nabla_{\q}u^\la]^2.\\
\endaligned\ee
Hence,
\be\nonumber\aligned
\mathcal{L}&:=L_1+L_2+L_3\\
&=-(n-4)\frac{d}{d\la}\int_{\pa B_1}\big[\la\frac{d}{d\la}(\nabla_{\q}u^\la)^2\big]-4(n-4)\frac{d}{d\la}\int_{\pa B_1}[\nabla_{\q}u^\la]^2\\
&+4(n-4)\la\int_{\pa B_1}\big[\frac{d}{d\la}\nabla_{\q}u^\la\big]^2\\
&=-(n-4)\frac{d}{d\la}\int_{\pa B_\la}\la^{6\frac{p+1}{p-1}-n-4}\frac{d}{d\la}\big[\la^2|\nabla u|^2-\la^2|\pa_r u|^2\big]\\
&\quad\quad-4(n-4)\frac{d}{d\la}\int_{\pa B_\la}\la^{6\frac{p+1}{p-1}-n-5}\big[\la^2|\nabla u|^2-\la^2|\pa_r u|^2\big]\\
&+4(n-4)\la\int_{\pa B_1}\big[\frac{d}{d\la}\nabla_{\q}u^\la\big]^2.\\
\endaligned\ee

\noindent{\bf Proof of Theorem \ref{10monoid}}. From the equation \eqref{10overlineE} and  combining with the estimates on $\mathcal{I},\mathcal{J},\mathcal{K},\mathcal{L}$ and the basic integrate by part, we obtain Theorem \ref{10monoid}.\hfill $\Box$

\section{Homogeneous solutions}

We first deduce polar coordinate representation for the  harmonic, biharmonic and triharmonic operator:

\be\nonumber
\Delta u=(\pa_{rr}+\frac{n-1}{r})u+\frac{1}{r^2}\Delta_\theta u,
\ee
and
\be\nonumber\aligned
\Delta^2u=&(\pa_{rr}+\frac{n-1}{r}\pa_{r})^2u+\Delta_\theta(\pa_{rr}+\frac{n-1}{r}\pa_r)(r^{-2}u)\\
&\;\;+\Delta_\theta(r^{-2}(\pa_{rr}+\frac{n-1}{r}\pa_r)u)+\Delta^2_\theta(r^{-4}u).
\endaligned\ee
Therefore,
\be\nonumber\aligned
\Delta^3u&=(\pa_{rr}+\frac{n-1}{r}\pa_r)^3u+\Delta_\theta^3(r^{-6}u)\\
&\quad+\Delta_\theta\Big((\pa_{rr}+\frac{n-1}{r}\pa_r)\big(r^{-2}(\pa_{rr}+\frac{n-1}{r}\pa_r)u\big)\\
&\quad+(\pa_{rr}+\frac{n-1}{r}\pa_r)^2(r^{-2}u)+r^{-2}(\pa_{rr}+\frac{n-1}{r}\pa_r)^2u\Big)\\
&\quad+\Delta_\theta^2\Big((\pa_{rr}+\frac{n-1}{r}\pa_r)(r^{-4}u)+r^{-4}(\pa_{rr}+\frac{n-1}{r}\pa_r)u\\
&\quad+r^{-2}(\pa_{rr}+\frac{n-1}{r}\pa_r)(r^{-2}u).
\endaligned\ee
Assume that $u$ is  homogeneous,  that is,  $u=r^{-\frac{6}{p-1}}w(\theta)$. By a direct calculation, the function  $w$ satisfy
\be\label{10ww}\aligned
\Delta_\theta^3w-k_2\Delta_\theta^2w+k_1\Delta_\theta w-k_0w=|w|^{p-1}w,
\endaligned\ee
where
\be\nonumber\aligned
&k_0=\\
&\frac{6}{p-1}(\frac{6}{p-1}+2)(\frac{6}{p-1}+4)(n-2-\frac{6}{p-1})(n-4-\frac{6}{p-1})(n-6-\frac{6}{p-1}),
\endaligned\ee
\be\nonumber\aligned
k_1=&\frac{6}{p-1}(\frac{6}{p-1}+2-n)(\frac{6}{p-1}+4)(\frac{6}{p-1}+6-n)\\
&\quad+(\frac{6}{p-1}+2)(\frac{6}{p-1}+4-n)(\frac{6}{p-1}+4)(\frac{6}{p-1}+6-n)\\
&\quad+\frac{6}{p-1}(\frac{6}{p-1}+2-n)(\frac{6}{p-1}+2)(\frac{6}{p-1}+4-n).
\endaligned\ee
\be\nonumber\aligned
k_2&=(\frac{6}{p-1}+4)(n-\frac{6}{p-1}-6)+\frac{6}{p-1}(n-\frac{6}{p-1}-2)\\
&\quad+(k+2)(n-\frac{6}{p-1}-4).
\endaligned\ee
Hence, test \eqref{10ww} with $w(\theta)$,  we have
\be\label{10II}\aligned
\int_{S^{n-1}}|\nabla_\theta\Delta_\theta w|^2+k_2|\Delta_\theta w|^2+k_1|\nabla_\theta w|^2+k_0w^2=\int_{S^{n-1}}|w|^{p+1}.
\endaligned\ee
For any $\varepsilon>0$, choose any $\eta_\varepsilon\in C_0^\infty(\frac{\varepsilon}{2},\frac{2}{\varepsilon})$, such that
$\eta_\varepsilon=1$ in $(\varepsilon,\frac{1}{\varepsilon})$ and
\be\nonumber
r|\eta'_\varepsilon(r)|+r^2|\eta''_\varepsilon(r)|+r^3|\eta'''_\varepsilon(r)|+r^4|\eta''''_\varepsilon(r)|+
r^5|\eta'''''_\varepsilon(r)|+
r^6|\eta''''''_\varepsilon(r)|\leq 1000.
\ee
Note that
\be\nonumber\aligned
\Delta(r^{-\frac{n-6}{2}}w(\theta)\eta_\varepsilon(r))=(\pa_{rr}+\frac{n-1}{r}\pa_r)(r^{-\frac{n-6}{2}}\eta_\varepsilon(r))w(\theta)
+r^{-\frac{n-2}{2}}\eta_\varepsilon(r)w(\theta)
\endaligned\ee
and that $|\nabla u|^2=\frac{1}{r^2}|\nabla_\theta u|^2+u_r^2$, we  see that
\be\nonumber\aligned
&|\nabla\Delta(r^{-\frac{n-6}{2}}w(\theta)\eta_\varepsilon(r))|^2\\
&=\frac{1}{r^2}|\nabla_\theta[(\pa_{rr}+\frac{n-1}{r}\pa_r)(r^{-\frac{n-6}{2}}\eta_\varepsilon(r))w(\theta)+r^{-\frac{n-2}{2}}\eta_\varepsilon(r)\Delta_\theta w]|^2\\
&\quad+|\frac{\pa}{\pa r}[(\pa_{rr}+\frac{n-1}{r}\pa r)(r^{-\frac{n-6}{2}}\eta_\varepsilon(r))w(\theta)+r^{-\frac{n-2}{2}}\eta_\varepsilon(r)\Delta_\theta w]|^2\\
&:=I_1+I_2.\\
\endaligned\ee
By a straightforward calculation we have the following estimates
\be\label{10I1}\aligned
\int_{\R^n}I_1\leq &\Big(\int_{S^{n-1}}|\nabla_\theta\Delta_\theta w|^2+\frac{(n-6)(n+2)}{2}|\Delta_\theta|^2+\frac{(n-6)^2(n+2)^2}{16}|\nabla_\theta w|^2d\theta\Big)\\
&\quad\cdot\Big(\int_0^\infty r^{-1}\eta_\varepsilon^2(r)dr\Big)\\
&+C\Big(\int_{S^{n-1}}|\nabla_\theta w|^2+|\Delta_\theta w|^2d\theta \Big)\\
&\quad\cdot\Big(\int_0^\infty \sum_{1\leq k+j\leq 4,k,j\geq0}r^{k+j-1}\eta^{(k)}_\varepsilon\eta^{(j)}_\varepsilon dr\Big)
\endaligned\ee
and
\be\label{10I2}\aligned
\int_{\R^n}I_2&\leq \Big( \int_{S^{n-1}}\frac{(n-2)^2}{4}|\Delta_\theta w|^2+\frac{(n-6)(n-2)^2(n+2)}{8}|\nabla_\theta w|^2\\
&\quad\quad+\frac{(n-6)^2(n-2)^2(n+2)^2}{64}w^2d\theta\Big) \cdot\Big( \int_0^\infty r^{-1}\eta^2_\varepsilon(r)dr\Big)\\
&+C\Big(\int_{S^{n-1}}w^2(\theta)+|\nabla_\theta w|^2\Big)\cdot \Big(\int_0^\infty \sum_{1\leq k+j\leq 6,k,j\geq0}r^{k+j-1}\eta^{(k)}_\varepsilon\eta^{(j)}_\varepsilon dr\Big).
\endaligned\ee
Here $\eta^{(k)}_\varepsilon:=\frac{d^k}{d r^k}\eta_\varepsilon$ for $k\geq1$ and $\eta^{(0)}_\varepsilon:=\eta_\varepsilon$.
Notice that

\be\nonumber\aligned
\int_0^\infty r^{-1}\eta^2_\varepsilon(r)dr\geq|\log\varepsilon|\rightarrow +\infty,\;\;\hbox{as}\;\; \varepsilon\rightarrow 0^+,
\endaligned\ee
\be\nonumber\aligned
\int_0^\infty \sum_{1\leq k+j\leq 6,k,j\geq0}r^{k+j-1}\eta^{(k)}_\varepsilon\eta^{(j)}_\varepsilon dr\leq C,
\endaligned\ee
where $C$ is independent of the radius $R$.
Recall the stability condition for triharmonic equation:
\be\nonumber\aligned
\int_{\R^n}|\nabla\Delta \phi|^2dx\geq p\int_{\R^n}|u|^{p-1}\phi^2dx.
\endaligned\ee
Let $\phi=r^{-\frac{n-6}{2}}w(\theta)\eta_\varepsilon(r)$. Then combining  with \eqref{10I1} and \eqref{10I2} and
letting $\varepsilon\rightarrow0$, we obtain that
\be\nonumber\aligned
p&\int_{s^{n-1}}|w|^{p+1}d\theta\leq\int_{s^{n-1}}|\nabla_\theta\Delta_\theta w|^2+\frac{3n^2-12n-20}{4}|\Delta_\theta w|^2\\
&+\frac{(n-6)(n+2)(3n^2-12n-4)}{16}|\nabla_\theta w|^2+\frac{(n-6)^2(n-2)^2(n+2)^2}{64}w^2(\theta),
\endaligned\ee
this combine with \eqref{10II}, we have the following estimate
\be\label{10w}\aligned
&\int_{s^{n-1}}(p-1)|\nabla_\theta\Delta_\theta w|^2+(pk_2-\frac{3n^2-12n-20}{4})|\Delta_\theta w|^2\\
&+(pk_1-\frac{(n-6)(n+2)(3n^2-12n-4)}{16})|\nabla_\theta w|^2\\
&+(pk_0-\frac{(n-6)^2(n-2)^2(n+2)^2}{64})w^2(\theta)\leq0.
\endaligned\ee
Since, $k=\frac{6}{p-1}$, we have $p=\frac{k+6}{k}$. Equivalently,  by the coefficient of the above inequality, we let
\be\label{10c}\aligned
c_0:&=(k+6)k_0/k-\frac{(n-6)^2(n-2)^2(n+2)^2}{64},\\
c_1:&=(k+6)k_1-\frac{(n-6)(n+2)(3n^2-12n-4)}{16}k,\\
c_2:&=(k+6)k_2-\frac{3n^2-12n-20}{4}k.
\endaligned
\ee
We consider the algebraic equation $c_0, c_1, c_2$ about the variable $k$, we only need consider positive real roots.
Roughly speaking, $c_0$ is a six-order algebraic equation about $k$, it has no explicit solution in general. Nonetheless we shall prove

\bl
Assume that $\frac{n+6}{n-6}<p<p_c (n)$. Then $ c_0, c_1, c_2 >0$.
\label{keyc123}
\el

Assuming the validity of Lemma \ref{keyc123}, we derive from (\ref{10w}) that $ w\equiv 0$. This gives
  \bt\label{10Thhomo} Let $u\in W_{loc}^{3,2}(\R^n\setminus \{0\})$ be a homogeneous, stable solution of \eqref{10LE}, for $\frac{n+6}{n-6}<p<p_c(n)$.
Assume that $|u|^{p+1}\in L_{loc}^1(\R^n\backslash\{0\})$, then $u \equiv 0.$
\et

 \vskip0.12in

The proof of Lemma \ref{keyc123} is quite  technical and thus we delay it to the appendix.

\section{Energy estimates and Blow down analysis }

 In the first part of this section, we obtain initial energy estimates on  the solutions of \eqref{10LE}, which are important when we perform a
 blow-down analysis in the second part of this  section.

 \subsection{Energy estimates}

\bl Let $u$ be a stable solution of \eqref{10LE}, then there exists a positive constant $C$ such that
\be\label{10estimate1}\aligned
&\int_{\R^n}|u|^{p+1}\eta^6+\int_{\R^{n}}|\nabla\Delta u|^2\eta^6\\
&\leq C\Big(\int_{\R^{n}}|\Delta u|^2\eta^4|\nabla\eta|^2+\int_{\R^{n}}|\nabla u|^2 \frac{|\Delta\eta^6|^2}{\eta^6}
+\int_{\R^{n}}u^2\frac{|\nabla\Delta \eta^6|^2}{\eta^6}\\
&\quad +\int_{\R^{n}}|\nabla u|^2\eta^2|\nabla \eta|^4
+\int_{\R^{n}}|\nabla^2 u|^2\eta^4|\nabla\eta|^2\\
&\quad +\int_{\R^{n}}|\nabla u|^2\eta^4|\nabla^2\eta|^2\Big).
\endaligned\ee
\el
\bp
Multiplying the equation \eqref{10LE} with $u\eta^6$, where $\eta$ is a test function, we get that
\be\label{10eee}\aligned
\int_{\R^n}|u|^{p+1}\eta^6&=\int_{\R^n}-\Delta^3 u\cdot u\eta^6=\int_{\R^n}\nabla\Delta^2u\cdot\nabla(u\eta^6)=-\int_{\R^n}\Delta^2u\Delta(u\eta^6)\\
&=\int_{\R^n}\nabla\Delta u\cdot\nabla\Delta(u\eta^6).
\endaligned\ee
Since $\Delta(\xi\eta)=\eta\Delta\xi+\xi\Delta\eta+2\nabla\xi\nabla\eta$, we have
\be\nonumber\aligned
\Delta(u\eta^6)=\eta^6\Delta u+u\Delta\eta^6+12\eta^5\nabla u\nabla\eta,
\endaligned\ee
therefore,
\be\label{10zhankai}\aligned
\nabla\Delta(u\eta^6)\nabla\Delta u=&6\eta^5\Delta u\nabla\eta\nabla\Delta u+(\eta)^6(\nabla\Delta u)^2+\Delta\eta^6\nabla u\nabla\Delta u\\
&+u\nabla\Delta \eta^6\nabla\Delta u+60\eta^4(\nabla\eta\nabla\Delta u)(\nabla u\nabla\eta)\\
&+12\eta^5\sum_{i,j}\pa_{ij}u\pa_i\eta\pa_j \Delta u+12\eta^5\sum_{i,j}\pa_i u\pa_{ij}\eta\pa_j\Delta u,
\endaligned\ee
where $\pa_j(j=1,...,n)$ denote the derivatives  with respect to  $x_1,...,x_n$ respectively.
A similar way can be applied to deal with the following term $|\nabla\Delta (u\eta^3)|^2$.
On the other hand, by the stability condition (see Definition \ref{zou-stable}), we have
\be\label{10stable1}
p\int_{\R^n}|u|^{p+1}\eta^6\leq\int_{\R^{n}}|\nabla\Delta(u\eta^3)|^2.
\ee
Combine  this with \eqref{10eee}, \eqref{10zhankai} and \eqref{10stable1}, we have
\be\nonumber\aligned
&\int_{\R^{n}}|\nabla\Delta u|^2\eta^6\\
&\leq  C\varepsilon\int_{\R^{n}}(\nabla\Delta u)^2\eta^6+C(\varepsilon)\big[\int_{\R^{n}}|\Delta u|^2\eta^4|\nabla\eta|^2\\
&\quad +\int_{\R^{n}}|\nabla u|^2 (\frac{|\Delta\eta^6|^2}{\eta^6}+\eta^4|\nabla^2\eta|^2)\\
&\quad +\int_{\R^{n}}u^2\frac{|\nabla\Delta \eta^6|^2}{\eta^6}
+\int_{\R^{n}}|\nabla u|^2\eta^2|\nabla \eta|^4
+\int_{\R^{n}}|\nabla^2 u|^2\eta^4|\nabla\eta|^2\big],\\
\endaligned\ee
we can select $\varepsilon$ so small  that $C\varepsilon\leq\frac{1}{2}$. Finally, combine with \eqref{10eee} and \eqref{10zhankai}, we obtain the conclusion of this lemma.
\ep


\bl\label{10LemmaE}
 Let $u$ be a stable solution of \eqref{10LE}.  Then

\be\label{10est1}\aligned
&\int_{B_{R}}|u|^{p+1}+\int_{B_{R}}(\nabla\Delta u)^2\leq C R^{-6}\int_{ B_{2R}} u^2,\endaligned\ee
\be\label{10est2}\aligned
&\int_{B_{R}}|u|^{p+1}+\int_{B_{R}}(\nabla\Delta u)^2\leq C R^{n-6\frac{p+1}{p-1}}.
\endaligned\ee
\el
\bp
We let $\eta=\xi^m$ where $m>1$ in the estimate \eqref{10estimate1}, we have
\be\label{10g}\aligned
&\int_{\R^n}|\nabla\Delta u|^2\xi^{6m}+\int_{\R^n}|u|^{p+1}\xi^{6m}\\
&\leq \int_{\R^n}u^2g_0(\xi)+\int_{\R^n}|\nabla u|^2 g_1(\xi)+\int_{\R^n}|\Delta u|^2g_2(\xi),
\endaligned\ee
where
\be\aligned
g_0(\xi):&=\xi^{6m-6}\sum_{0\leq i+j+k+r+s+t\leq 6}|\nabla^i \xi||\nabla^j \xi||\nabla^k \xi||\nabla^r \xi||\nabla^s \xi||\nabla^t \xi|,\\
g_1(\xi):&=\xi^{6m-4}\sum_{0\leq i+j+k+l\leq 4}|\nabla^i \xi||\nabla^j \xi||\nabla^k \xi||\nabla^l \xi|,\\
g_2(\xi):&=\xi^{6m-2}\sum_{0\leq i+j\leq 2}|\nabla^i \xi||\nabla^j \xi|,
\endaligned\ee
where we define $\nabla^0\xi:=\xi$ and notice that $g_m(\xi)\geq0$ for $m=0,1,2$. Now, we claim that
\be\aligned
g_1^2(\xi)\leq C g_0(\xi)g_2(\xi),\;\;|\nabla^2g_2(\xi)|\leq C g_1(\xi),\;\;g_2^2(\xi)\leq C\xi^{6m}g_1(\xi).
\endaligned\ee
This claim can be verified  by direct calculations and  will be used for the following estimates.
Since $|\nabla u|^2=\frac{1}{2}\Delta(u^2)-u\Delta u$, we have
\be\label{10g1}\aligned
\int_{\R^n}|\nabla u|^2g_1(\xi)&=\frac{1}{2}\int_{\R^n}\Delta(u^2)g_1(\xi)-\int_{\R^n}u\Delta u g_1(\xi)\\
&=\frac{1}{2}\int_{\R^n}u^2\Delta g_1(\xi)-\int_{\R^n}u\Delta u g_1(\xi)\\
&\leq \frac{1}{2}\int_{\R^n}u^2\Delta g_1(\xi)+\varepsilon\int_{\R^n}(\Delta u)^2g_2(\xi)+\frac{1}{4\varepsilon}\int_{\R^n}u^2 g_0(\xi).
\endaligned\ee
We note the following  differential identity
\be\nonumber
(\Delta u)^2=\sum_{j,k}(u_j u_k)_{jk}-\sum_{j,k}(u_{jk})^2-2\nabla\Delta u\cdot\nabla u.
\ee
Hence
$(\Delta u)^2\leq\sum_{j,k}(u_j u_k)_{jk}-2\nabla\Delta u\cdot\nabla u.$
Therefore we have
\be\label{10g2}\aligned
\int_{\R^n}(\Delta u)^2g_2(\xi)&\leq\int_{\R^n}\sum_{j,k}(u_j u_k)_{jk}g_2(\xi)-2\int_{\R^n}\nabla\Delta u\cdot\nabla g_2(\xi)\\
&=\int_{\R^n}\sum_{j,k}u_j u_k g_2(\xi)_{jk}-2\int_{\R^n}\nabla\Delta u\cdot\nabla g_2(\xi)\\
&\leq C\int_{\R^n}|\nabla u|^2g_1(\xi)+\delta\int_{\R^n}|\nabla\Delta u|^2\xi^{6m}+C(\delta)\int_{\R^n}|\nabla u|^2g_1(\xi)\\
&\leq C\int_{\R^n}|\nabla u|^2g_1(\xi)+\delta\int_{\R^n}|\nabla\Delta u|^2\xi^{6m}.
\endaligned\ee
Combining with \eqref{10g1} and \eqref{10g2}, by selecting  the positive parameter $\varepsilon$ small enough, we can obtain that
\be\nonumber\aligned
\int_{\R^n}|\nabla u|^2g_1(\xi)+\int_{\R^n}(\Delta u)^2g_2(\xi)\leq C\int_{\R^n}u^2g_0(\xi)+\delta\int_{\R^n}|\nabla\Delta u|^2\xi^{6m}.
\endaligned\ee
By    combining the above inequalities with \eqref{10g} and selecting  the positive parameter $\delta$ small enough, we have  that
\be
\int_{\R^n}|\nabla\Delta u|^2\xi^{6m}+\int_{\R^n}|u|^{p+1}\xi^{6m}\leq C\int_{\R^n}u^2g_0(\xi).
\ee
This proves \eqref{10est1}. Further, we let $\xi=1$ in $B_{R}$ and $\xi=0$ in $B_{2R}^C$, satisfying $|\nabla \xi |\leq\frac{C}{R}$, we have
\be\aligned
&\int_{\R^n}|\nabla\Delta u|^2\xi^{6m}+\int_{\R^n}|u|^{p+1}\xi^{6m}\leq C\int_{\R^n}u^2g_0(\xi)\leq CR^{-6}\int_{\R^n}u^2\xi^{6m-6}\\
&\leq CR^{-6}(\int_{\R^n}|u|^{p+1}\xi^{(3m-3)(p+1)})^{\frac{2}{p+1}}R^{n(1-\frac{2}{p+1})}.
\endaligned\ee
By selecting $m>1$ and letting $m$ close to 1, we can make sure that $(3m-3)(p+1)\leq 6m$. It follows  that   \eqref{10est2} holds.
\ep

\vskip0.1in

\subsection{ Blow-down analysis and  the proof of Theorem \ref{10Liouville}}

\vskip0.13in
  {\bf The proof of Theorem \ref{10Liouville}.}   Firstly, we consider    $1<p\leq\frac{n+6}{n-6}$. If $p<\frac{n+6}{n-6}$, we can let $R\rightarrow \infty$ in \eqref{10est2} to get $u\equiv0$ directly. However,  if $p=\frac{n+6}{n-6}$, this gives
\be\nonumber
\int_{\R^n}|\nabla\Delta u|^2+|u|^{p+1}<+\infty.
\ee
Hence
\be\nonumber
\lim_{R\rightarrow+\infty}\int_{B_{2R}(x)\setminus B_R(x)}|\nabla\Delta u|^2+|u|^{p+1}=0.
\ee
Then by \eqref{10est2}, and noting that now $n=6\frac{p+1}{p-1}$, we have
\be\nonumber\aligned
&\int_{B_R(x)}|\nabla\Delta u|^2+|u|^{p+1}\leq CR^{-6}\int_{B_{2R}(x)\setminus B_R(x)}u^2\\
&\leq C R^{-6}(\int_{B_{2R}(x)\setminus B_R(x)}|u|^{p+1})^{\frac{2}{p+1}}R^{n(1-\frac{2}{p+1})}\leq C(\int_{B_{2R}(x)\setminus B_R(x)}|u|^{p+1})^{\frac{2}{p+1}}.
\endaligned\ee
Let $R\rightarrow+\infty$, we get that $u\equiv0$.

\vskip0.15in

Secondly, we consider the supercritical case, i.e., $p>\frac{n+6}{n-6}$. We complete the proof via a few steps.

\vskip0.1in

\noindent{\bf Step 1.}  $\lim_{\la\rightarrow\infty} E(u,0,\la)<\infty$.\\

From Theorem \ref{10Monotone} we know that $E$ is nondecreasing w. r. t. $\la$, so we only need to show that $E(u,0,\la)$ is bounded. Note that
\be\nonumber
E(u,0,\la)\leq\frac{1}{\la}\int_{\la}^{2\la} E(u,0,t)dt\leq\frac{1}{\la^2}\int_\la^{2\la}\int_{t}^{t+\la}E(u,0,\gamma)d\gamma dt.
\ee
From Lemma \ref{10LemmaE}, we have that
\be\nonumber\aligned
\frac{1}{\la^2}&\int_\la^{2\la}\int_{t}^{t+\la}\gamma^{6\frac{p+1}{p-1}-n}\big[\int_{B_\gamma}\frac{1}{2}|\nabla\Delta u|^2dx -\frac{1}{p+1}\int_{B_\gamma}|u|^{p+1}dx\big]d\gamma dt\leq C,
\endaligned\ee
where $C>0$ is independent of $\gamma$.
\be\label{10type1}\aligned
\frac{1}{\la^2}&\int_\la^{2\la}\int_{t}^{t+\la}\int_{\pa B_\gamma}\gamma^{6\frac{p+1}{p-1}-n-5}\big[\frac{6}{p-1}(\frac{6}{p-1}-1)(\frac{6}{p-1}-2)u\\
&+\frac{18}{p-1}(\frac{6}{p-1}-1)\gamma\pa_r u
+\frac{18}{p-1}\gamma^2 \pa_{rr} u+\gamma^3\pa_{rrr} u\big]\\
&\quad\quad\big[\frac{6}{p-1}(\frac{6}{p-1}-1) u+\frac{12}{p-1}\gamma\pa_r u+\gamma^2\pa_{rr}u\big]\\
&\leq C\frac{1}{\la^2}\int_\la^{2\la}\int_{t}^{t+\la}t^{6\frac{p+1}{p-1}-n-5}\int_{\pa B_\gamma} \big[u^2+\gamma^2(\pa_r u)^2+\gamma^4(\pa_{rr}u)^2+\gamma^6(\pa_{rrr}u)^2\big]\\
&\leq C\frac{1}{\la^2}\int_\la^{2\la}t^{6\frac{p+1}{p-1}-n-5}\int_{B_{3\la}}\big[u^2+\gamma^2(\pa_r u)^2+\gamma^4(\pa_{rr}u)^2+\gamma^6(\pa_{rrr}u)^2\big]\\
&\leq C \la^{n-6\frac{p+1}{p-1}+6}\frac{1}{\la^2}\int_\la^{2\la}t^{6\frac{p+1}{p-1}-n-5}dt\\
&\leq C
\endaligned\ee
and

\be\label{10type2}\aligned
&\Big|\frac{1}{\la^2}\int_\la^{2\la}\int_{t}^{t+\la}\int_{\pa B_\gamma}\gamma^{6\frac{p+1}{p-1}-n-4}\frac{d}{d\gamma}(\gamma^2\Delta u-\gamma^2\pa_{rr}u-(n-1)\gamma\pa_r u)^2\Big|\\
&\leq\frac{1}{\la^2}\int_\la^{2\la}t^{6\frac{p+1}{p-1}-n-5}\int_{t}^{t+\la}\int_{\pa B_\gamma}|\big[2\gamma^2\Delta u-2\gamma^2\pa_{rr}u-(n-1)\gamma\pa_r u\big]\\
&\quad\big[\gamma^2\Delta u-\gamma^2\pa_{rr}u-(n-1)\gamma\pa_r u\big]|\\
&\leq\frac{1}{\la^2}\int_\la^{2\la}t^{6\frac{p+1}{p-1}-n-5}\int_{B_{3\la}}|[2\gamma^2\Delta u-2\gamma^2\pa_{rr}u-(n-1)\gamma\pa_r u]\\
&\quad \big[\gamma^2\Delta u-\gamma^2\pa_{rr}u-(n-1)\gamma\pa_r u\big]|\\
&\leq C\frac{1}{\la^2}\int_\la^{2\la}t^{6\frac{p+1}{p-1}-n-5}\int_{B_{3\la}}
\big[\gamma^4(\Delta u)^2+\gamma^4(\pa_{rr}u)^2+\gamma^2\pa_r u\big]\\
&\leq C \la^{n-6\frac{p+1}{p-1}+6}\frac{1}{\la^2}\int_\la^{2\la}t^{6\frac{p+1}{p-1}-n-5}dt\\
&\leq C.
\endaligned\ee
The remaining    terms can be treated similarly as the estimate \eqref{10type1} or \eqref{10type2}.

\vskip0.32in

\noindent{\bf Step 2.} For any $\la>0$, recall the definition
\be\nonumber
u^\la(x):=\la^{\frac{6}{p-1}}u(\la x),
\ee
and $u^\la$ is also a smooth solution of \eqref{10LE} on $\R^n$. By rescaling \eqref{10est2}, for $\la>0$ and balls $ B_r(x)\subset\R^n$,
\be\nonumber
\int_{B_r(x)}|\nabla\Delta u^\la|^2+|u^\la|^{p+1}\leq C r^{n-6\frac{p+1}{p-1}}.
\ee
In particular, $u^\la$ are uniformly bounded in $L^{p+1}_{loc}(\R^n)$ and $\nabla\Delta u^\la$ are uniformly bounded in $L^2_{loc}(\R^n)$.
By elliptic estimates, $u^\la$ are also uniformly bounded in $W^{3,2}_{loc}(\R^n)$. Hence up to a subsequence of $\la\rightarrow+\infty$,
we can assume that $u^\la\rightarrow u^\infty$ weakly in $W^{3,2}_{loc}(\R^n)\cap L^{p+1}_{loc}(\R^n)$. By Sobolev embedding, $u^\la\rightarrow u^\infty$
strongly in $W^{2,2}_{loc}(\R^n)$. Then for any ball $B_R(0)$, by interpolation and noting \eqref{10est2}, for any $q\in [1,p+1)$ as $\la\rightarrow+\infty$,
\be\label{10LL}
\|u^\la-u^\infty\|_{L^q(B_R(0))}\leq\|u^\la-u^\infty\|_{L^1(B_R(0))}^t\|u^\la-u^\infty\|_{L^{p+1}(B_R(0))}^{1-t}\rightarrow0,
\ee
where $t\in (0,1]$ satisfies  $\frac{1}{q}=t+\frac{1-t}{p+1}$. That is, $u^\la\rightarrow u^\infty$ in $L^q_{loc}(\R^n)$ for any $q\in [1,p+1)$.
For any function $\phi\in C_0^\infty(\R^n)$, we have that
\be\nonumber
\int_{\R^n}\nabla\Delta u^\infty\cdot \nabla \Delta \phi-(u^\infty)^{p-1}u^\infty \phi=\lim_{\la\rightarrow+\infty}\int_{\R^n}\nabla\Delta u^\la\cdot \nabla \Delta \phi-(u^\la)^{p-1}u^\infty \phi,
\ee
\be\nonumber
\int_{\R^n}|\nabla\Delta \phi|^2-p(u^\infty)^{p-1}\phi^2=\lim_{\la\rightarrow+\infty}\int_{\R^n}|\nabla\Delta \phi|^2-p(u^\la)^{p-1}\phi^2.
\ee
Therefore $u^\infty\in W^{3,2}_{loc}(\R^n)\cap L_{loc}^{p+1}(\R^n)$ is a stable solution of \eqref{10LE} in $\R^n$.
\vskip0.12in

\noindent{\bf Step 3.}  {\bf The function $u^\infty$ is homogeneous.}
Due to the scaling invariance of the functional $E $  (i.e., $E(u,0,R\la)=E(u^{\la},0,R)$ )
and the monotonicity formula, for any given $R_2>R_1>0$, we see that
\be\nonumber\aligned
0=&\lim_{i\rightarrow\infty}\big(E(u,0,R_2\la_i)-E(u,0,R_1\la_i)\big)\\
=&\lim_{i\rightarrow\infty}\big(E(u^{\la_i},0,R_2)-E(u^{\la_i},0,R_1)\big)\\
\geq&C(n,p)\liminf_{i\rightarrow\infty}\int_{B_{R_2}\setminus B_{R_1}}
 r^{6\frac{p+1}{p-1}-n-6}\big(\frac{6}{p-1}{u^{\la_i}}+r\frac{\pa u^{\la_i}}{\pa r}\big)^2dydx\\
 \geq&C(n,p)\int_{B_{R_2}\setminus B_{R_1}}
 r^{6\frac{p+1}{p-1}-n-6}\big(\frac{6}{p-1}{u^{\infty}}+r\frac{\pa u^{\infty}}{\pa r}\big)^2dydx.\\
\endaligned\ee
In the last inequality we have used the weak convergence of  the sequence $(u^{\la_i})$ to the function $u^{\infty}$ in $W^{1,2}_{loc}(\R^n)$
as $i\to \infty$. This implies that
\be\nonumber
\frac{6}{p-1}\frac{u^{\infty}}{r}+\frac{\pa u^{\infty}}{\pa r}=0\;\;\hbox{a.e.}\;\;\hbox{in}\;\;\R^{n}.
\ee
Integrating in $r$ shows that
\be\nonumber
u^\infty(x)=|x|^{-\frac{6}{p-1}}u^\infty(\frac{x}{|x|}).
\ee
That is,  $u^\infty$ is homogeneous.

\vskip0.12in

\noindent{\bf Step 4. $u^\infty=0$.}
This is a direct consequence of Theorem \ref{10Thhomo}.
Since this holds for the limit of any sequence $\la\rightarrow+\infty$, by \eqref{10LL} we get
\be\nonumber
\lim_{\la\rightarrow+\infty} u^\la\quad \hbox{strongly in}\;\; L^2(B_4(0)).
\ee

\vskip0.12in

\noindent{\bf Step 5. $u=0$.}  For all $\la\rightarrow+\infty$, we see that
 \be\nonumber
 \lim_{\la\rightarrow+\infty}\int_{B_4(0)}(u^\la)^2=0.
 \ee
By \eqref{10est2},
\be\label{10s1}\lim_{\la\rightarrow+\infty}\int_{B_3(0)}|\nabla\Delta u^\la|^2+|u^\la|^{p+1}\leq\lim_{\la\rightarrow+\infty}\int_{B_4(0)}(u^\la)^2=0.\ee
By the elliptic interior $L^2$ estimate, we get that
\be\nonumber
\lim_{\la\rightarrow+\infty}\int_{B_2(0)}\sum_{k\leq 3}|\nabla^k u^\la|^2=0.
\ee
In particular, we can choose a sequence $\la_i\rightarrow+\infty$ such that
\be\nonumber
\int_{B_2(0)}\sum_{k\leq 3}|\nabla^k u^{\la_i}|^2\leq 2^{-i}.
\ee
Hence we have
\be\nonumber
\int_1^2\sum_{i=1}^{+\infty}\int_{\pa B_r}\sum_{k\leq 3}|\nabla^k u^{\la_i}|^2dr\leq\sum_{i=1}^{+\infty}\int_1^2\int_{\pa B_r}\sum_{k\leq 3}|\nabla^k u^{\la_i}|^2dr\leq1.
\ee
That is, the function
\be\nonumber
H(r):=\sum_{i=1}^{+\infty}\int_{\pa B_r}\sum_{k\leq3}|\nabla^k u^{\la_i}|^2\in L^1(1,2).
\ee
Then there exists an $r_0\in (1,2)$ such that $H(r_0)<+\infty$, by which  we get that
\be\nonumber
\lim_{i\rightarrow+\infty}\|u^{\la_i}\|_{W^{3,2}(\pa B_{r_0})}=0.
\ee
Combining this with \eqref{10s1} and the scaling invariance of $E(r)$, we have
\be\nonumber
\lim_{i\rightarrow+\infty}E(\la_i r_0, 0, u)=\lim_{i\rightarrow+\infty}E(r_0, 0, u^{\la_i})=0.
\ee
Since $\la_i r_0\rightarrow+\infty$ and $E(r, 0, u)$ is non-decreasing in $r$, we get
\be\nonumber
\lim_{r\rightarrow+\infty}E(r, 0, u)=0.
\ee
By the smoothness of $u$, $\lim_{r\rightarrow0}E(r, 0,u)=0$. Then again by the monotonicity of $E(r, 0, u)$ and step 4, we obtain that
\be\nonumber
E(r, 0, u)=0\quad\hbox{for all}\quad r>0.
\ee
Therefore,  by the monotonicity formula we know that  $u$ is homogeneous, then $u\equiv0$ by Theorem \ref{10Thhomo}. \hfill$\Box$

\section{Finite Morse index solutions}
In this section, we prove Theorem \ref{10Liouvillec}. We always assume that $u$ is a smooth classical solution
with finite Morse index.

\bl\label{10LemmaF1}
Let $u$ be  a smooth  solution (positive or sign changing) of \eqref{10LE} with finite Morse index, then there exist   constants $C>0$  and $R_0$ such that
\be\nonumber
|u(x)|\leq C|x|^{-\frac{6}{p-1}}, \quad  \forall x\in B_{R_0}(0)^c.
\ee

\el
\bp
Since  that $u$ is stable outside $B_{R_0}$. For $x\in B_{R_0}^c$, let $M(x)=|u(x)|^{\frac{p-1}{6}}$ and $d(x)=|x|-R_0$. Assume that there exists a sequence of
$x_k\in B_{R_0}^c$ such that
\be\label{10M}
M(x_k)d(x_k)\geq2k.
\ee
Since $u$ is bounded on any compact set of $\R^n$, $d(x_k)\rightarrow+\infty$.

By the doubling Lemma \cite{Souplet2007}, there exists another sequence $y_k\in B_{R_0}^c$, such that
\be\nonumber\aligned
&M(y_k)d(y_k)\geq2k,\quad M(y_k)\geq M(x_k);\\
&M(z)\leq2M(y_k)\;\;\hbox{for any}\;\;z\in B_{R_0}^c \;\;\hbox{such that}\;\;|z-y_k|\leq\frac{k}{M(y_k)}.
\endaligned\ee
Now we define
\be\nonumber
u_k(x):=M(y_k)^{-\frac{6}{p-1}}u(y_k+M(y_k)^{-1}x),\quad\hbox{for}\;x\in B_k(0).
\ee
This and above arguments  give that, $u_k(0)=1$, $|u_k|\leq 2^{\frac{6}{p-1}}$ in $B_k(0)$.
 Further,  $B_{k/M(y_k)}\cap B_{R_0}=\emptyset$,
which implies that $u$ is a stable solution in $B_{k/M(y_k)}(y_k)$. Hence, $u_k$ is stable in $B_k(0)$.

By elliptic regularity theory, $u_k$ are uniformly bounded in $C^7_{loc}(B_k(0))$. Up to a sequence,
$u_k$ convergence to $u_\infty$ in $C^6_{loc}(\R^n)$. By the above conditions on $u_k$, we have
\be\aligned
|u_{\infty}(0)|=1,\;\; |u_\infty|\leq 2^{\frac{6}{p-1}};\;\; u_\infty \hbox{is a smooth stable solution of \eqref{10LE} in }\;\;\R^n.
\endaligned\ee
By the Liouville theorem for stable solution, we have $u_\infty\equiv0$, a contradiction with \eqref{10M}.
\ep

\bc\label{10CF1}
There exist   constants $C>0$ and $R_0$ such that for all $x\in B_{R_0}^c$,
\be\label{10estd}
\sum_{k\leq5}|x|^{\frac{6}{p-1}+k}|\nabla^ku(x)|\leq C.
\ee
\ec

\bp
For any $x_0$ with $|x_0|>3R_0$, take $\la=\frac{|x_0|}{2}$ and define
\be\nonumber
\overline{u}(x)=\la^{\frac{6}{p-1}}u(x_0+\la x).
\ee
By the previous Lemma, $\overline{u}\leq C$ in $B_1(0)$. By the elliptic regularity theory we have
\be\nonumber
\sum_{k\leq5}|\nabla^k \overline{u}(0)|\leq C.
\ee
Scaling back we get the conclusions.
\ep

\subsection{The proof of Theorem \ref{10Liouvillec}-(1):  $1<p<\frac{n+6}{n-6}$ (Subcritical case)}

We need the following Pohozaev identity. A general version can be seen in \cite{Wei1999}.

\bl
For any function $u$ satisfying \eqref{10LE}, we have that
\be\nonumber
(\frac{n-6}{2}-\frac{n}{p+1})\int_{B_R}|u|^{p+1}=\int_{\pa B_R}B_3(u)d\sigma,
\ee
where
\be\aligned
-B_3(u)=&\frac{R}{p+1}|u|^{p+1}-2(-\Delta)^2u\frac{\pa u}{\pa n}+2u\frac{\pa (-\Delta)^2u}{\pa n}-\frac{R}{2}|\nabla\Delta u|^2\\
&-\frac{n-2}{2}u\frac{\pa\Delta^2 u}{\pa n}+\frac{n-2}{2}\Delta u\frac{\pa \Delta u}{\pa n}+<x\cdot\nabla u>\frac{\pa \Delta^2 u}{\pa n}\\
&-\Delta^2u\frac{\pa<x\cdot\nabla u>}{\pa n}+<x\cdot\nabla\Delta u>\frac{\pa \Delta u}{\pa n}.
\endaligned\ee
\el

\noindent{\bf The proof of Theorem \ref{10Liouvillec}-(1).} By Corollary \ref{10CF1}, for $R>R_0$ ($R_0$ is  defined in Corollary \ref{10CF1}), noting that $p<\frac{n+6}{n-6}$ (hence $n-6\frac{p+1}{p-1}<0$), we have the following estimate
\be\nonumber
\int_{\pa B_R}|B_3(u)|d\sigma\leq C\int_{\pa B_R}R^{-\frac{12}{p-1}-5}d\sigma\leq CR^{n-6\frac{p+1}{p-1}}\rightarrow0\;\hbox{as}\;R\rightarrow+\infty.
\ee
Combining with the Pohazaev identity, letting $R\rightarrow+\infty$, we get that
\be\nonumber
(\frac{n-6}{2}-\frac{n}{p+1})\int_{\R^n}|u|^{p+1}=0.
\ee
Since $\frac{n-6}{2}-\frac{n}{p+1}<0$, we obtain that $u\equiv0$.

\subsection{The proof of Theorem \ref{10Liouvillec}-(3):  $p=\frac{n+6}{n-6}$ (critical case)}

Since $u$ is stable outside $B_{R_0}$, Lemma \ref{10LemmaE} still holds if the support of $\eta$ is outside $B_{R_0}$.
Take $\varphi\in C_0^\infty(B_{2R}\setminus B_{2R_0})$ such that $\varphi\equiv1$ in $B_{R}\setminus B_{3R_0}$ and
$\sum_{k\leq5}|x|^k|\nabla^k u|\leq 1000$. Then by choosing $\eta=\varphi^m$, where $m$ is bigger than 1, we get that
\be\nonumber
\int_{B_R\setminus B_{3R_0}}|\nabla\Delta u|^2+|u|^{p+1}\leq C.
\ee
Letting $R\rightarrow+\infty$, we have
\be\label{10LF1}
\int_{\R^n}|\nabla\Delta u|^2+|u|^{p+1}<+\infty.
\ee
By the interior elliptic estimates and Holder's inequality, we have

\be\nonumber\aligned
R^{-4}\int_{B_{2R}\setminus B_R}|\nabla u|^2&\leq C\int_{B_{3R}\setminus B_{R/2}}|\nabla\Delta u|^2+C\Big(\int_{B_{3R}\setminus B_{R/2}}|u|^{p+1}\Big)^{\frac{2}{p+1}},\\
R^{-2}\int_{B_{2R}\setminus B_R}|\Delta u|^2&\leq C\int_{B_{3R}\setminus B_{R/2}}|\nabla\Delta u|^2+C\Big(\int_{B_{3R}\setminus B_{R/2}}|u|^{p+1}\Big)^{\frac{2}{p+1}},\\
R^{-6}\int_{B_{2R}\setminus B_R}u^2&\leq C\int_{B_{3R}\setminus B_{R/2}}|\nabla\Delta u|^2+C\Big(\int_{B_{3R}\setminus B_{R/2}}|u|^{p+1}\Big)^{\frac{2}{p+1}}.
\endaligned\ee
Therefore, we have that
\be\nonumber
\max\Big(R^{-4}\int_{B_{2R}\setminus B_R}|\nabla u|^2,R^{-2}\int_{B_{2R}\setminus B_R}|\Delta u|^2,R^{-6}\int_{B_{2R}\setminus B_R}u^2\Big)\rightarrow0
\ee
 as  $R\rightarrow+\infty.$ On the other hand, testing \eqref{10LE} with an compact support function $\eta^2$, we get
\be\nonumber
\int_{\R^n}|\nabla\Delta u|^2\eta^2-|u|^{p+1}\eta^2=-\int_{\R^n}\nabla\Delta u\cdot\nabla\Delta\eta^2\cdot u+\nabla\Delta u\nabla u\Delta\eta^2
+\nabla\Delta u\nabla(2\nabla u\nabla\eta^2).
\ee
By selecting $\eta(x)=\xi(\frac{x}{R})^{3m}$, $m>1$ and $\xi\in C_0^\infty(B_{2})$ and $\xi\equiv1$ in $B_1$, and $\sum_{k\leq3}|\nabla^k u|\leq1000$, we get that
\be\nonumber\aligned
&\Big|\int_{\R^n}|\nabla\Delta u|^2\xi(\frac{x}{R})^{6m}-|u|^{p+1}\xi(\frac{x}{R})^{6m}\Big|\leq C\Big(R^{-4}\int_{B_{2R}\setminus B_R}|\nabla u|^2\\
&\quad\quad\quad +R^{-2}\int_{B_{2R}\setminus B_R}|\Delta u|^2+R^{-6}\int_{B_{2R}\setminus B_R}u^2\Big).
\endaligned\ee
Now letting $R\rightarrow+\infty$, we obtain that
\be\nonumber
\int_{\R^n}|\nabla\Delta u|^2-|u|^{p+1}=0.
\ee
Combining with \eqref{10LF1}, we get the conclusions.

\subsection{The proof of Theorem \ref{10Liouvillec}-(2):   $p>\frac{n+6}{n-6}$ (supercritical case)}

\bl
There exists a  constant $C>0$ such that $E(r, 0,u)\leq C$  for all $r>3R_0$.
\el

\bp
From the monotonicity formula, combine the derivative estimates \eqref{10estd}, we have the following estimates
\be\nonumber\aligned
E(r, 0,u)\leq& Cr^{4\frac{n+6}{n-6}-n}\Big(\int_{B_r}|\nabla\Delta u|^2+|u|^{p+1}\Big)\\
&+\sum_{j,k\leq4,j+k\leq5}r^{6\frac{p+1}{p-1}-n-5+j+k}\int_{\pa B_r}|\nabla^j u||\nabla^k u|\\
\leq &C.
\endaligned\ee
This constant only depends on the constant in \eqref{10estd}.
\ep
As a consequence, we have the following
\bc
\be\label{cor-zou-1}
\int_{B_{3R_0}^c}\frac{(\frac{6}{p-1}u(x)+|x|u_r(x))^2}{|x|^{n+6-6\frac{n+6}{n-6}}}dx<+\infty.
\ee
\ec
As before, we define the blowing down sequence
\be\nonumber
u^\la(x)=\la^{\frac{6}{p-1}}u(\la x).
\ee
By Lemma \ref{10LemmaF1}, $u^\la$ are uniformly bounded in $C^7(B_r(0)\setminus B_{1/r}(0))$ for any fixed $r>1$ and moreover,
$u^\la$ is stable outside $B_{R_0/\la}(0)$. There exists  a function $u^\infty\in C^6(\R^n\setminus\{0\})$,
such that up to a subsequence of $\la\rightarrow+\infty$, $u^\la$ converges to $u^\infty$ in $C^6_{loc}(\R^n\setminus\{0\})$,
$u^\infty$ is stable solution of \eqref{10LE} in $\R^n\setminus\{0\}$.

 For any $r>1$, by the above Corollary \ref{cor-zou-1},
\be\nonumber\aligned
 &\int_{B_r\setminus B_{1/r}}\frac{(\frac{6}{p-1}u^\infty(x)+|x|u^\infty_r(x))^2}{|x|^{n+6-6\frac{n+6}{n-6}}}dx\\
 &=\lim_{\la\rightarrow+\infty}\int_{B_r\setminus B_{1/r}}\frac{(\frac{6}{p-1}u^\la(x)+|x|u^\la_r(x))^2}{|x|^{n+6-6\frac{n+6}{n-6}}}dx\\
 &=\lim_{\la\rightarrow+\infty}\int_{B_{\la r}\setminus B_{\la/r}}\frac{(\frac{6}{p-1}u(x)+|x|u_r(x))^2}{|x|^{n+6-6\frac{n+6}{n-6}}}dx\\
 &=0.
\endaligned\ee
Hence, $u^\infty$ is homogeneous, and by Theorem \ref{10Thhomo}, $u^\infty\equiv0$ if $p<p_c(n)$.
 Since this hold for any limit of $u^\la$ as $\la\rightarrow+\infty$, then we have
 \be\nonumber
 \lim_{|x|\rightarrow+\infty}|x|^{\frac{6}{p-1}}|u(x)|=0.
 \ee
Then as in the proof of Corollary \ref{10CF1}, we have
\be\nonumber
\lim_{|x|\rightarrow+\infty}\sum_{k\leq 6}|x|^{\frac{6}{p-1}+k}|\nabla^k u(x)|=0.
\ee
Therefore,  for any $\varepsilon>0$, take an $R_0$ such that for $|x|>R_0$,  there holds
\be\nonumber
\sum_{k\leq 6}|x|^{\frac{6}{p-1}+k}|\nabla^k u(x)|\leq \varepsilon.
\ee
Then for any $r\gg R_0$, we have
\be\nonumber\aligned
&E(r;0,u)\leq C r^{6\frac{p+1}{p-1}-n}\int_{B_R(0)}\big(|\nabla\Delta u|^2+|u|^{p+1}\big)\\
&+C\varepsilon r^{6\frac{p+1}{p-1}-n}\int_{B_r(0)\setminus B_R(0)}|x|^{-6\frac{p+1}{p-1}}
+C\varepsilon r^{6\frac{p+1}{p-1}+1-n}\int_{\pa B_r(0)}|x|^{-6\frac{p+1}{p-1}}\\
&\leq C(R_0)(r^{6\frac{p+1}{p-1}-n}+\varepsilon).
\endaligned\ee
We obtain that $\lim_{r\rightarrow+\infty}E(r;0,u)=0$ since $6\frac{p+1}{p-1}-n<0$ and $\varepsilon$ can be arbitrarily small.
On the other hand,  since $u$ is smooth we have $\lim_{r\rightarrow0}E(r;0,u)=0$, thus $E(r;0,u)=0$ for all $r>0$, thus by the monotonicity formula we get that
$u$ is homogeneous, and then by Theorem \ref{10Thhomo}, we know  that $u\equiv0$.
\section{Proofs of Theorem \ref{10Monotone}}

To further investigate the optimal condition to make the monotonicity formula hold, we find that
we have drop the term $\int_{\pa B_1}\la(\frac{dv^\la}{d\la})^2$. Recall that

\be\nonumber\aligned
v^\la=\Delta u^\la&=\la^2\frac{d^2u^\la}{d\la^2}+(n-1-\frac{12}{p-1})\la\frac{du^\la}{d\la}+\frac{6}{p-1}(\frac{6}{p-1}-n+2)u^\la+\Delta_\theta u^\la\\
&:=\la^2\frac{d^2u^\la}{d\la^2}+a\la\frac{du^\la}{d\la}+bu^\la+\Delta_\theta u^\la.
\endaligned\ee
By some integrate by part, we have that

\be\nonumber\aligned
&\int_{\pa B_1}\la(\frac{dv^\la}{d\la})^2=\int_{\pa B_1}\Big(\la^5(\frac{d^3u^\la}{d\la^3})^2+(a^2-2a-2b-4)\la^3(\frac{d^2u^\la}{d\la^2})^2\\
&\quad\quad\quad\quad\quad\quad\quad\quad+(-a^2+b^2+2a+2b)\la(\frac{du^\la}{d\la})^2\Big)\\
&+\int_{\pa B_1}\Big(-2\la^3(\nabla_\theta \frac{d^2u^\la}{d\la^2})^2+(10-2b)\la(\nabla_\theta \frac{du^\la}{d\la})^2\Big)
+\int_{\pa B_1}\la(\Delta_\theta \frac{du^\la}{d\la})^2\\
&+\frac{d}{d\la}\Big(\int_{\pa B_1}\sum_{0\leq i,j\leq 2,i+j\leq2}{c^1_{i,j}}\la^{i+j}\frac{d^iu^\la}{d\la^i}\frac{d^ju^\la}{d\la^j}
+\sum_{0\leq s,t\leq 2,s+t\leq2}{c^2_{s,t}}\la^{s+t}\frac{d^su^\la}{d\la^s}\frac{d^tu^\la}{d\la^t}\Big)
\endaligned\ee
where $c^i_{i,j},c^2_{s,t}$ determined by $a,b$ hence by $p,n$.
From Theorem \ref{10monoid} we derive that
\be\label{10Ec}\aligned
&\frac{d E^c(\la,x,u)}{d\la}=\int_{\pa B_1}3\la^5\Big(\frac{d^3u^\la}{d\la^3}\Big)^2+A_{1}\la^3\Big(\frac{d^2u^\la}{d\la^2}\Big)^2+A_{2}\la\Big(\frac{du^\la}{d\la}\Big)^2\\
&+\int_{\pa B_1}\Big(2\la^3(\nabla_\theta \frac{d^2u^\la}{d\la^2})^2+(8\alpha-4\beta-2b+4n-18)\la(\nabla_\theta \frac{du^\la}{d\la})^2\Big)\\
&+\int_{\pa B_1}\la(\Delta_\theta \frac{du^\la}{d\la})^2,\\
\endaligned\ee
where
\be\label{10A1A2}\aligned
A_1:&=10\delta_1-2\delta_2-56+a^2-2a-2b-4\\
A_2:&=-18\delta_1+6\delta_2-4\delta_3+2\delta_4+72-a^2+b^2+2a+2b,
\endaligned\ee
Let $k:=\frac{6}{p-1}$, a direct calculation we have that
\be\nonumber\aligned
A_1=-10k^2+(-60+10n)k-n^2+24n-83,
\endaligned\ee
\be\nonumber\aligned
A_2&=3k^4+(36-6n)k^3+(3n^2-48n+150)k^2\\
&\quad+(12n^2-114n+252)k+9n^2-72n+135,
\endaligned\ee
and
\be\nonumber\aligned
B_1:=8\alpha-4\beta-2b+4n-18=-6k^2+(-36+6n)k+12n-42.
\endaligned\ee

Notice that our supercritical condition $p>\frac{n+6}{n-6}$ is equivalent to $0<k<\frac{n-6}{2}$.

Firstly, we have the following lemma which yields the sign of $A_2$ and $B_1$.

\bl\label{10A2-0} If $p>\frac{n+6}{n-6}$, then $A_2>0$ and $B_1>0$.
\el

\bp From  \eqref{10A1A2}, we derive that
\be\label{10A2}
A_2=3(k+1)(k+3)(k-(n-5))(k-(n-3)),
\ee
and the roots of $B_1=0$ are

\be\nonumber
\frac{1}{2}n-3-\frac{1}{2}\sqrt{n^2-4n+8},\frac{1}{2}n-3+\frac{1}{2}\sqrt{n^2-4n+8}
\ee
Recall that $p>\frac{n+6}{n-6}$ is equivalent to $0<k<\frac{n-6}{2}$, we get the conclusion.
\ep

 To show monotonicity formula,  we  proceed it to prove the following inequality
\be
\label{keyinequality}
3\la^5(\frac{d^3u^\la}{d\la^3})^2+A_{1}\la^3\Big(\frac{d^2u^\la}{d\la^2}\Big)^2+A_{2}\la(\frac{du^\la}{d\la})^2\ee
\be
\geq \epsilon\la(\frac{du^\la}{d\la})^2+\frac{d}{d\la}
\Big(\sum_{0\leq i,j\leq2}c_{i,j}\la^{i+j}\frac{d^iu^\la}{d\la^i}\frac{d^ju^\la}{d\la^j}
\Big).
\ee
To deal with the rest of the dimensions, we employ the second idea: we find nonnegative constants $d_1, d_2$ and constants $c_1, c_2$ such that we have the following Jordan form decomposition:
\be\label{10quad}\aligned
&3\la^5(f''')^2+A_1\la^3(f'')^2+A_2\la(f')^2=3\la(\la^2f'''+c_1\la f'')^2+d_1\la(\la f''+c_2 f')^2\\
&\;+d_2\la(f')^2
+\frac{d}{d\la}(\sum_{i,j}e_{i,j}\la^{i+j}f^{(i)}f^{(j)}),
\endaligned\ee
where the  unknown constants are to be determined.

\bl
\label{60011}
Let $ p>\frac{n+6}{n-6}$ and $A_1$ satisfy
\be
\label{A11}
 A_1+12 >0,
 \ee
 then there exist  nonnegative numbers $d_1,d_2$,    and real numbers $c_1,c_2,e_{i,j}$ such that
the  differential inequality (\ref{10quad}) holds.

 \el

 \bp Since $$4\la^4f'''f''=\frac{d}{d\la}(2\la^4(f'')^2)-8\la^3(f'')^2$$ and
 $$2\la^2f''f'=\frac{d}{d\la}(\la^2(f')^2)-2\la(f')^2,$$ by comparing the coefficients of $\lambda^3(f'')^2$ and $\lambda (f')^2,$  we have that

\be\nonumber
d_1=A_1-3c_1^2+12c_1, \quad d_2=A_2-(c_2^2-2c_2)(A_1-3c_1^2+12c_1).
\ee
In particular, $$\max_{c_1}{d_1(c_1)}=A_1+12   \hbox{  and the critical point is }   c_1=2. $$
Since $A_2>0$,  we select that $c_1=2,c_2=0$.  Hence,  in this case,
by a direct calculation we see that $d_1=A_1+12>0$.  Then we get  the conclusion.
\ep

We conclude from Lemma \ref{60011} that if $A_1+12>0$ then  (\ref{keyinequality}) holds.  This implies that when $7\leq n\leq20, p>\frac{n+6}{n-6}$ or $n\geq 21$ and
\be
\label{pee}
\frac{n+6}{n-6}<p<\frac{5n+30-\sqrt{15n^2-60n+190}}{5n-30-\sqrt{15n^2-60n+190}}
\ee
 then (\ref{keyinequality}) holds.
\vskip0.2in

Combine idea from the above with the the following idea, we can get better condition to make the monotonicity formula holds.
We start from the differential identity \eqref{10quad}. Recall that  the derivative term is a 'good' term since it can be absorbed in the term $E^c(\la,x,u)$.

We make use of two ideas to prove (\ref{keyinequality}).  The second idea is straightforward. We   use  the  positivity  of terms $A_2\la(\frac{du}{d\la})^2$ and $3\la^5(\frac{d^3u}{d\la^3})^2$ to bound the
  term $A_1(\frac{d^2u}{d\la^2})^2$. Note
\be\nonumber\aligned
3&\la(\la^2f'''+2\la f'')^2+A_2\la(f')^2\geq-2\sqrt{3A_2}(\la^3f'''f'+2\la^2 f''f')\\
&=2\sqrt{3A_2}\Big(\la^3(f'')^2-\la(f')^2)\Big)+2\sqrt{3 A_2}\frac{d}{d\la}\Big(\la^3 f''f'-\frac{1}{2}\la^2(f')^2\Big).
\endaligned\ee

 We observe that
We divide the term $A_2\la(f')^2$ into two parts, $\theta A_2(f')^2$
and $(1-\theta)\la(f')^2$, and then,  find the optimal parameter $\theta$. Following this idea, we have
\be\aligned
&3\la(\la^2f'''+2\la f'')^2+(A_1+12)\la^3(f'')^2+A_2\la(f')^2\\
&\geq (A_1+12+2\sqrt{3\alpha A_2})\la^3(f'')^2+\Big((1-\alpha)A_2-2\sqrt{3\alpha A_2}\Big)\la(f')^2\\
&\quad +2\sqrt{3\alpha A_2}\frac{d}{d\la}\Big(\la^3 f''f'-\frac{1}{2}\la^2(f')^2\Big),\\
\endaligned\ee
hence, we have the desired    monotonicity formula  once the following two inequalities hold:
\be\label{10mean2}
A_1+12+2\sqrt{3\alpha A_2}\geq0, \;\; (1-\alpha)A_2-2\sqrt{3\alpha A_2}>0.
\ee
The second inequality of \eqref{10mean2} gives the range of $\alpha$, that is
\be\label{10alpha}
\frac{12\alpha}{(\alpha-1)^2}<\min_{0\leq k\leq\frac{n-6}{2}}A_2.
\ee
Obviously, the first inequality of \eqref{10mean2} holds if $A_1+12\geq0$. So we just need to check the following inequality:
\be
\label{A13}
 A_1+12 >-\sqrt{12 \alpha A_2} \ \mbox{where} \  \ \frac{8\alpha}{(\alpha-1)^2}<\min_{0\leq k\leq\frac{n-6}{2}}A_2.
\ee

We discuss the remaining dimensions as follows:

When $n=21$,$$\min_{0\leq k\leq\frac{n-6}{2}\mid n=14}A_2=A_2(k=0)=2592,$$ thus from \eqref{10alpha} we get that $\alpha\leq0.9342$, then $(A_1+12)^2-12\alpha\mid_{\alpha=0.9342}A_2<0$ if $-0.5941782055<k< 4.483334837$. On the other hand, $A_1+12<0$ implies that $0<k<0.05352432355$.
Hence, (\ref{A13}) holds.

The case of $ 21\leq n\leq 30$ can be dealt with similarly. We omit the details here.

\medskip

Let
\be\nonumber
p_m(n):=\begin{cases}
+\infty\;\;\;\;\;\;\;\;&\hbox{if}\;\;\;\;\;\;\;\; n\leq30,\\
\frac{5n+30-\sqrt{15n^2-60n+190}}{5n-30-\sqrt{15n^2-60n+190}}&\hbox{if}\;\;\;\;\;\;\;\;n\geq31.\\
\end{cases}
\ee
Combining  all the lemmas of this section,  we obtain   the following theorem.
\bt\label{LZ=62609} For $\frac{n+6}{n-6}<p<p_m(n)$,  then there exists  a $ C(n,p)>0$ such that
\be\nonumber
\frac{d}{d\la}E^c(\la,x,u)\geq C(n,p)\int_{\pa B_1}\la(\frac{du^\la}{d\la})^2.
\ee
\et

\vskip0.2in

\noindent{\bf Proof of  Theorem \ref{10Monotone}.}  Let $d(n)$ be defined at (\ref{dndef})(See also the appendix).  By Lemma \ref{10h2} of the appendix, we know that $d(n)<\sqrt{n}$ for $n\geq15$. Hence, we have the following inequalities
\be\label{10nt-0}
\frac{n-8}{2}-d(n)>\frac{n-8}{2}-\sqrt{n}
 \ee
and
\be\label{10nt}
\frac{n-8}{2}-\sqrt{n}\geq\frac{1}{2}n-3-\frac{1}{10}\sqrt{15n^2-60n+190}\;\;\hbox{for}\;\;n\geq14.\ee
Hence we derive that  when  $n\geq 15$
\be\nonumber
\frac{n+4-2d(n)}{n-8-2d(n)}<\frac{5n+30-\sqrt{15n^2-60n+190}}{5n-30-\sqrt{15n^2-60n+190}}.
\ee
Therefore, $p_c(n)< p_m(n)$. Theorem \ref{10Monotone} is thus proved.  \hfill$\Box$

\section{Appendix: Proof of Lemma \ref{keyc123}}

In this appendix, we prove the technical lemma \ref{keyc123}.

Recall the definition of $c_0, c_1$ and $c_2$ in (\ref{10c}). Let
 $$k:=\frac{n-8}{2}+a.$$
Then $c_0$ can be rewritten in terms of $a$:
\be\nonumber\aligned
c_0=&-a^6+(8+\frac{3}{4}n^2)a^4-(16+\frac{3}{16}n^4)a^2\\
&\;\;+\frac{3}{16}n^5-\frac{15}{16}n^4-\frac{3}{2}n^3+\frac{33}{4}n^2+3n-9.
\endaligned
\ee
Further, we let $t:=a^2$, we get a three-order algebraic equation as following:
\be\label{10c0t}\aligned
c_0=&-t^3+(8+\frac{3}{4}n^2)t^2-(16+\frac{3}{16}n^4)t\\
&\;\;+\frac{3}{16}n^5-\frac{15}{16}n^4-\frac{3}{2}n^3+\frac{33}{4}n^2+3n-9.
\endaligned
\ee
By the two crucial transformation above, we reduce a six-order algebraic equation to a third order algebraic equation.
Now we can get the explicit solution of the above equation \eqref{10c0t} which  has two imaginary   roots and one real root.  We denote  the real root as $d(n)$.  Let
\be\nonumber\aligned
d_1(n):=-94976+20736n+103104n^2-10368n^3+1296n^5-3024n^4-108n^6,
\endaligned\ee
\be\nonumber\aligned
d_2(n):&=6131712-16644096n^2+6915840n^4-690432n^6-3039232n\\
&\quad+4818944n^3-1936384n^5+251136n^7-30864n^8-4320n^9\\
&\quad+1800n^{10}-216n^{11}+9n^{12}
\endaligned\ee
and
\be\label{10d01}
d_0(n):=-(d_1(n)+36\sqrt{d_2(n)})^{1/3}.
\ee
Notice that
\be\nonumber\aligned
d_2(n):=&(9n^8-216n^7+1872n^6-6048n^5-16032n^4+206208n^3\\
&\;\;-848640n^2-189952n+383232)(n-2)^2(n+2)^2>0 \;\;\hbox{if}\;\;n\geq12,
\endaligned\ee
hence $\sqrt{d_2(n)}$ is well defined  whenever $n\geq12$.  Define
\be\label{10dn}
d(n):=\frac{1}{6}\Big(9n^2+96-\frac{1536+1152n^2}{d_0(n)}-\frac{3}{2}d_0(n)\Big)^{1/2}.
\ee
By the proof of Lemma \ref{10h2} below, we will see that   $d(n)$ is well-defined,  i.e., $9n^2+96>\frac{1536+1152n^2}{d_0(n)}+\frac{3}{2}d_0(n)$.

\vskip0.23in
Let $r_1,r_2$ denote the two real roots of $c_0$ which can be computed as
\be\label{10r}
r_1:=\frac{n-8}{2}-d(n),\quad r_2:=\frac{n-8}{2}+d(n).
\ee
Therefore,   we see that  $c_0>0$  whenever  $r_1<k<r_2$.
Since the  roots $r_1$ and $r_2$ depend on  $d_0(n)$,   we must have a fine estimate on $d_0(n).$
\bl\label{10h1}  The $d_0(n)$ has the following properties:
\begin{itemize}
\item [(1)]
\be\label{10d02}
d_0(n):=\frac{256(3n^2+4)}{(36\sqrt{d_2(n)}-d_1(n))^{1/3}}.
\ee

\item [(2)]  For $n\geq15$, then
\be\nonumber \frac{d}{d n}d_0(n)<0,\;\; 128<d_0(n)<187.\ee
\end{itemize}

\el
\bp The proof of $(1)$ comes from the following identity which can be checked directly:
\be\label{10L68}
d_1^2(n)-36^2d_2^2(n)=256^3(3n^2+4)^3.
\ee
Now we start  to prove the conclusion  $(2)$ of the Lemma.
From \eqref{10L68} we see that $d_0(n)>0$. Thus we only need to  show that $\frac{d}{dn}d^3_0(n)<0$. In fact,
\be\label{10Lddn}\aligned
\frac{d}{d n}d^3_0(n)=-\frac{d}{d n}d_1(n)-\frac{18\frac{d}{d n}d_2(n)}{\sqrt{d_2(n)}},
\endaligned\ee
since
\be\nonumber
-\frac{d}{d n}d_1(n)=648n^5-6480n^4+12096n^3+31104n^2-206208n-20736>0\;\;\hbox{if}\;\;n\geq8,
\ee
and furthermore
\be\nonumber\aligned
&\Big(-\frac{d}{d n}d_1(n)\Big)^2\cdot d_2(n)-18^2\Big(\frac{d}{d n}d_2(n)\Big)^2\\
=&-293534171136n^{15}+4109478395904n^{14}-9001714581504n^{13}\\
&-168292924784640n^{12}+1233104438034432n^{11}-3119550711201792n^{10}\\
&-6748415824232448n^9+21348066225291264n^8
-1783991975804928n^7\\
&+9835612546793472n^6+34945090870837248n^5-114643053771227136n^4\\
&+19014404334944256n^3
-110880250103070720n^2-14427791579676672n\\
&-356241767399424\\
=&-10871635968(3n^3-18n^2+84n+8))(n^4-8n^3-40n^2+480n+16)\\
&\cdot (n-2)^2(n+2)^2(3n^2+4)^2<0\;\;\hbox{if}\;\;n\geq3.
\endaligned\ee
Thus,  combining the above two equations,  we get that
\be\nonumber
-\frac{d}{d n}d_1(n)\cdot \sqrt{d_2(n)}<18\frac{d}{d n}d_2(n).
\ee
Hence, if we combine this with \eqref{10Lddn}, we have that
\be\nonumber\aligned
\frac{d}{d n}d^3_0(n)=-\frac{d}{d n}d_1(n)-\frac{18\frac{d}{d n}d_2(n)}{\sqrt{d_2(n)}}<0.
\endaligned\ee
Therefore,  $\frac{d}{dn}d_0(n)<0$ for $n\geq15$. Notice that $d_0(n)\mid_{n=15}\simeq 186.0929<187$
and a straightforward calculation shows that
\be\nonumber
\lim_{n\rightarrow+\infty}d_0(n)=128.
\ee
By the monotonicity of $d_0(n)$ for $n\geq15$, we derive that $d_0(n)\in(128,187)$ for $n\geq15$.

\ep

By straightforward calculation we have the following asymptotic properties.
\bl\label{10h3}
\be\nonumber\lim_{n\rightarrow+\infty}\frac{d(n)}{\sqrt{n}}=1,\quad\lim_{n\rightarrow+\infty}(d(n)-\sqrt{n})=0,\ee
\be\nonumber\lim_{n\rightarrow+\infty}\sqrt{n}(d(n)-\sqrt{n})=-\frac{1}{2}.\ee
\el
By  Lemma  \ref{10h3} above, we known that $d(n)$ behaves as $\sqrt{n}-\frac{1}{2}\frac{1}{\sqrt{n}}$ if  $n$ large. Although $\lim_{n\rightarrow+\infty}\frac{d(n)}{\sqrt{n}}=1$,
 the limit   behavior  gives no information on the  size relation  between $d(n)$ and $\sqrt{n}$. Therefore, we  need the following more  delicate analysis.

\bl\label{10h2} When  $n\geq15$, we have
\be\nonumber
d(n)<\sqrt{n}.
\ee
\el

\bp We prove the following inequality: For $n\geq15$,

\be\label{10Ld0n}
9n^2-36n+96<\frac{1536+1152n^2}{d_0(n)}+\frac{3}{2}d_0(n)<9n^2+96.
\ee
The second  inequality of \eqref{10Ld0n} is  equivalent to the following
\be\label{10Ld0n1}
x^2-(6n^2+64)x+768n^2+1024<0.
\ee
Here $x\in (128,187)$ since $d_0(n)\in(128,187)$. Next, we show that \eqref{10Ld0n1} holds.
The roots of the equation  corresponding to the above inequality  are
\be\nonumber
x_1(n)=3n^2+32-3n\sqrt{n^2-64},\;\;\; x_2(n)=3n^2+32+3n\sqrt{n^2-64}.
\ee
For $n\geq15$, we have that $x_2(n)\geq x_2(15)\geq1276>d_0(n)$. Next we will show that
$d_0(n)>x_1(n)$, which implies that \eqref{10Ld0n1} holds, hence  the second  inequality  of \eqref{10Ld0n} holds.
Since $n^2-64>(n-3)^2$ for $n\geq15$, we have that
\be\nonumber\aligned
0<x_1(n)=&3n^2+32-3n\sqrt{n^2-64}=\frac{768n^2+1024}{3n^2+32+3n\sqrt{n^2-64}}\\
<&\frac{768n^2+1024}{3n^2+32+3n(n-3)}=\frac{768n^2+1024}{6n^2-9n+32}.
\endaligned\ee
To compare $x_1(n)$ with $d_0(n)$, let us compare $x^3_1(n)$ with $d^3_0(n)$.
First, we know  that

\be\label{10x11}\aligned
&-d_1(n)(6n^2-9n+32)^3-(768n^2+1024)^3\\
&=23328n^{12}-384912n^{11}+2443608n^{10}-8266860n^9-276048n^8+76177584n^7\\
&\quad-915397632n^6+1095581376n^5-4004833536n^4+1592960256n^3\\
&\quad-2731991040n^2-3305373696n+2038431744>0\;\;\hbox{if}\;\;n\geq10.
\endaligned\ee
It follows  that

\be\label{10x12}\aligned
&\Big(-d_1(n)(6n^2-9n+32)^3-(768n^2+1024)^3\Big)^2-36^2d_2(n)(6n^2-9n+32)^6\\
&=1358954496\Big(116640n^{17}-606528n^{16}+1195560n^{15}+16771860n^{14}\\
&\quad-104564844n^{13}+682366923n^{12}
-1464330096n^{11}+5142941100n^{10}\\
&\quad-6506609472n^9+15562840464n^8-11332244736n^7
+21360207936n^6\\
&\quad-5590593536n^5+10574331904n^4+4294279168n^3-2878341120n^2\\
&\quad+3791650816n-3221225472\Big)\\
&=1358954496(4320n^{11}-22464n^{10}+27000n^9+711036n^8-4003812n^7\\
&\quad+22548513n^6-38373440n^5+96546304n^4-66202112n^3+68272128n^2\\
&\quad+59244544n-50331648))(3n^2+4)^3>0\;\;\hbox{if}\;\;n\geq1.
\endaligned\ee
Then combining with \eqref{10x11} and \eqref{10x12}, we get that
\be\nonumber\aligned
-d_1(n)(6n^2-9n+32)^3-(768n^2+1024)^3>36\sqrt{d_2(n)}(6n^2-9n+32)^3.
\endaligned\ee
Hence,

\be\nonumber
-d_1(n)-36\sqrt{d_2(n)}>\frac{(768n^2+1024)^3}{(6n^2-9n+32)^3},
\ee
that is $d^3_0(n)>x_1^3(n)$, which yields that  $d_0(n)>x_1(n)$.  Combining with $d_0(n)<x_2(n)$  when  $n\geq15$, we obtain that \eqref{10Ld0n1}. Hence, we get the
second inequality of \eqref{10Ld0n}.
\vskip0.225in

A similar technique can be applied to the first inequality of   \eqref{10Ld0n}, which is
 equivalent  to the following inequality:
\be\label{10Ld0n2}
x^2-(6n^2-24n+64)x+768n^2+1024>0,
\ee
where $x\in (128,187)$ since $d_0(n)\in(128,187)$.
The roots of the equation  corresponding to the above inequality  are
\be\nonumber\aligned
r_1(n)&=3n^2-12n+32-\sqrt{9n^4-72n^3-432n^2-768n};\\
r_2(n)&=3n^2-12n+32+\sqrt{9n^4-72n^3-432n^2-768n}.
\endaligned\ee
Notice that $9n^4-72n^3-432n^2-768n>0$ if $n\geq13$.
We next shows that $d_0(n)<r_1(n)<r_2(n)$, hence we get \eqref{10Ld0n2}.
For $r_1(n)$, notice that
\be\nonumber\aligned
9n^4-72n^3-432n^2-768n&<9n^4-72n^3-432n^2+2304n+9216\\
&=(3n^2-12n-96)^2,
\endaligned\ee
it follows that

\be\nonumber\aligned
r_1(n)=&\frac{768n^2+1024}{3n^2-12n+32+\sqrt{9n^4-72n^3-432n^2-768n}}\\
>&\frac{768n^2+1024}{3n^2-12n+32+3n^2-12n-96}=\frac{384n^2+512}{3n^2-12n-32}:=r_{10}(n).
\endaligned\ee
Notice that $3n^2-12n-32=3(n+4)(n-8)>0$ if $n\geq9$.

Firstly  we observe  that
\be\label{10r11}\aligned
&-d_1(n)(3n^2-12n-32)^3-(384n^2+512)^3\\
&=2916n^{12}-69984n^{11}+548208n^{10}-699840n^9-12052800n^8+54991872n^7\\
&\;\;-7831296n^6-691006464n^5-299151360n^4+4048994304n^3+3403284480n^2\\
&\;\;-2821718016n-3246391296>0\;\;\hbox{if}\;\;n\geq11,
\endaligned\ee
and that
\be\label{10r12}\aligned
&\Big(-d_1(n)(3n^2-12n-32)^3-(384n^2+512)^3\Big)^2-36^2d_2(n)(3n^2-12n-32)^6\\
&=5435817984\Big(-729n^{16}+14580n^{15}-36936n^{14}-631152n^{13}+3184272n^{12}\\
&\;\;+6849792n^{11}
-15453504n^{10}-49876992n^9-32256000n^8-28111872n^7\\
&\;\;+268692480n^6+613150720n^5+898416640n^4+1187315712n^3+983040000n^2\\
&\;\;+616562688n+369098752
\Big)\\
&=-5435817984(27n^{10}-540n^9+1260n^8+25536n^7-123120n^6-352960n^5\\
&\;\;+1058048n^4+3124224n^3-2383872n^2-9633792n-5767168)(3n^2+4)^3\\
&\;\;<0\;\;\hbox{if}\;\;n\geq14.
\endaligned\ee
Combining with \eqref{10r11} and \eqref{10r12}, we have that

\be\nonumber\aligned
-d_1(n)(3n^2-12n-32)^3-(384n^2+512)^3<36\sqrt{d_2(n)}(3n^2-12n-32)^3,
\endaligned\ee
hence

\be\nonumber\aligned
-d_1(n)-36\sqrt{d_2(n)}<\frac{(384n^2+512)^3}{(3n^2-12n-32)^3},
\endaligned\ee
that is, $d^3_0(n)<r^3_{10}(n)$, thus $d_0(n)<r_{10}(n)<r_1(n)<r_2(n)$.  Therefore, \eqref{10Ld0n2} holds.  This is, the
first inequality of  \eqref{10Ld0n} holds.

\vskip0.105in

Summing up,  recall that \eqref{10dn}, $d(n)$ is well defined and  in particular, $ d(n)<\sqrt{n}$ for $n\geq15$.
These are direct consequences of \eqref{10Ld0n}.

\ep

\vskip0.12in
Next,  we will show that, under the case $\frac{n+6}{n-6}<p<p_c(n)$, we have $c_0>0,c_1>0,c_2>0$ simultaneously.
Notice that $\frac{n+6}{n-6}<p<p_c(n)$ equivalent to $\min\{0,r_1(n)\}<k<\frac{n-6}{2}$. And $r_1(n)$ is exactly the root of $c_0=0$, hence $c_0>0$ directly. The condition $\min\{0,r_1(n)\}<k<\frac{n-6}{2}$ is not very applicable, in view of the estimates in Lemma \ref{10h2} and Lemma \ref{10h3} above, we can evaluate  $c_1,c_1$ under the interval $\frac{n-8}{2}-\sqrt{n}<k<\frac{n-8}{2}+\sqrt{n}$. In the two lemmas below, we will  follow this idea.

\bl\label{10h4} Under the condition $\frac{n-8}{2}-\sqrt{n}<k<\frac{n-8}{2}+\sqrt{n}$,  for $n\geq36$, we have
\be\nonumber c_1>0.
\ee
\el
\bp

Combining with \eqref{10c}, we have
\be\nonumber\aligned
c_{1}&=3k^5+(54-6n)k^4+(3n^2-84n+372)k^3+(30n^2-408n+1224)k^2\\
&\quad +(\frac{159}{2}n^2-810n+1917-\frac{3}{16}n^4+\frac{3}{2}n^3)k+48n^2-480n+1152.
\endaligned\ee
Set $k=\frac{n-8}{2}+a(n)\sqrt{n}$,  where $-1\leq a(n)\leq1$. For  the simplicity, we denote   $a(n)$ by $a$. Thus,
\be\nonumber\aligned
c_{1}&=12+(\frac{9}{8}-\frac{3}{4}a^2)n^4+(-\frac{3}{2}a^3+\frac{3}{2}a)n^{\frac{7}{2}}+(-\frac{39}{4}+\frac{3}{2}a^4+3a^2)n^3\\
&\quad+(3a^5-\frac{3}{2}a)n^{\frac{5}{2}}
+(-6a^4+6a^2+3)n^2+(-12a^3-18a)n^{\frac{3}{2}}\\
&\quad+(24a^2+\frac{141}{2})n-3an^{\frac{1}{2}}.
\endaligned\ee
For the case $0\leq a\leq1$, since $3a^5-\frac{3}{2}a\geq-\frac{3}{10^{\frac{5}{4}}}$, we get from the above identity that
\be\nonumber\aligned
c_{1}&\geq 12+\frac{3}{8}n^4-\frac{39}{4}n^3-\frac{3}{10^{\frac{5}{4}}}n^{\frac{5}{2}}-3n^2-30n^{\frac{3}{2}}+\frac{141}{2}n-3n^{\frac{1}{2}}\\
&=\frac{3}{8}t^8-\frac{39}{4}t^6-\frac{3}{10^{\frac{5}{4}}}-3t^4-30t^3+\frac{141}{2}t^2-3t+12\quad (n=t^2)\\
&\geq0\quad\hbox{if}\quad n=t^2\geq 26.8(t\geq5.168).
\endaligned\ee
For the case $-1\leq a\leq0$, since $\frac{3}{2}(a-a^3)\geq-\frac{\sqrt{3}}{3}$, we have
\be\nonumber\aligned
c_{1}&\geq12+\frac{3}{8}n^4-\frac{\sqrt{3}}{3}n^{\frac{7}{2}}-\frac{39}{4}n^3-\frac{3}{2}n^{\frac{5}{2}}-3n^2-30n^{\frac{3}{2}}+\frac{141}{2}n\\
&=\frac{3}{8}t^8-\frac{\sqrt{3}}{3}t^7-\frac{39}{4}t^6-\frac{3}{2}t^5-3t^4-30t^3+\frac{141}{2}+12\\
&\geq0\quad\hbox{if}\quad n=t^2\geq 35.98(t\geq5.999).
\endaligned\ee

\ep

\bl\label{10h5} Under the condition $\frac{n-8}{2}-\sqrt{n}<k<\frac{n-8}{2}+\sqrt{n}$, for $n\geq12$, we have
\be\nonumber c_2>0.
\ee
\el
\bp
We set $k=\frac{n-8}{2}+a(n)\sqrt{n}$, hence by the assumption we have $-1\leq a(n)\leq1$. For simplicity, we denote $a$ for $a(n)$.
From \eqref{10c}, we have
\be\aligned
c_2:=-3k^2+(-36+3n)k^2+(-135+27n-\frac{3}{4}n^2)k+36n-192.
\endaligned\ee
Plugging $k=\frac{n-8}{2}+a\sqrt{n}$ in the above expression we get that
\be\nonumber\aligned
c_{2}=-36+(\frac{9}{2}-\frac{3}{2}a^2)n^2+(-3a^3+3a)n^{\frac{3}{2}}-\frac{39}{2}n+9an^{\frac{1}{2}}.
\endaligned\ee
If  $0\leq a\leq1$, we have
\be\nonumber\aligned
c_{2}&\geq -36+3n^2-3n^{\frac{3}{2}}-\frac{39}{2}n\\
&=3t^4-3t^3-\frac{39}{2}t^2-36\quad (n=t^2)\\
&\geq 0 \quad\hbox{if}\quad n=t^2\geq 10.9025(t\geq3.3019).
\endaligned\ee
If  $-1\leq a\leq 0$, we have
\be\aligned
c_{2a}&\geq -36+3n^2-3n^{\frac{3}{2}}-\frac{39}{2}n-9n^{\frac{1}{2}}\\
&=3t^4-3t^3-\frac{39}{2}t^2-9t-36\quad (n=t^2)\\
&\geq 0 \quad\hbox{if}\quad n=t^2\geq 11.8259(t\geq3.4388).
\endaligned\ee
\ep
Next,  we state a lemma  via the   numerical analysis of the above arguments.
\bl\label{10h6}

Consider the supercritical case $p>\frac{n+6}{n-6}$, i.e.,  $0<k<\frac{n-6}{2}$. We have the following facts.
\begin{itemize}
\item [(1)] If  $ 0<k<\frac{n-6}{2}$ and $n\leq 14$, then $c_0,c_1,c_2>0$;
\item [(2)] If  $ 15\leq n\leq 50$ and $r_1<k<\frac{n-6}{2}$, then $c_0,c_1,c_2>0$.

\end{itemize}
\el

\vskip0.3in

Notice that $k>\min\{r_1:=\frac{n-8}{2}-d(n),0\}$ is equivalent to $ p<p_c (n)$. Combining Lemmas \ref{10h4}-\ref{10h6} we obtain the proof of Lemma \ref{keyc123}.

\bigskip


\end{document}